\documentclass[10pt]{article}

\RequirePackage{amsmath}
\RequirePackage{amsfonts}
\RequirePackage{amssymb}
\RequirePackage{hyperref}
\RequirePackage{multirow}
\RequirePackage{algorithm}
\RequirePackage{algpseudocode}
\RequirePackage{graphicx}
\RequirePackage{xcolor}
\RequirePackage{amsthm}
\RequirePackage{amsmath}
\RequirePackage{amsfonts}
\RequirePackage{amssymb}
\RequirePackage{hyperref}
\RequirePackage{multirow}
\RequirePackage{algorithm}
\RequirePackage{algpseudocode}
\RequirePackage{graphicx}
\RequirePackage{xcolor}
\allowdisplaybreaks

\usepackage{dpr_fec,ifthen,subcaption,dsfont}

\usepackage{dsfont}

\usepackage{isetechreport}
\def\reportyear{20T}
\def\reportno{009}
\def\revisionno{1}
\def\originaldate{May 15, 2020}
\def\revisiondate{\today}
\coraltrue
\cvcrfalse

\newcommand{\citep}{\cite}
\newcommand{\Halmos}{\qed}
\newcommand{\argmin}{\arg\min}
\newcommand{\argmax}{\arg\max}
\newenvironment{henumerate}{\benumerate}{\eenumerate}
\newenvironment{hitemize}{\bitemize}{\eitemize}

\begin{document}

\title{Gradient Sampling Methods with Inexact Subproblem Solves and Gradient Aggregation\thanks{This material is based upon work supported by the National Science Foundation under grant numbers CCF--1618717 and CCF--1740796.}}

\author{Frank E.~Curtis\thanks{E-mail: frank.e.curtis@lehigh.edu}}
\author{Minhan Li\thanks{E-mail: mil417@lehigh.edu}}
\affil{Department of Industrial and Systems Engineering, Lehigh University}
\titlepage

\maketitle

\begin{abstract}
  Gradient sampling (GS) has proved to be an effective methodology for the minimization of objective functions that may be nonconvex and/or nonsmooth.  The most computationally expensive component of a contemporary GS method is the need to solve a convex quadratic subproblem in each iteration.  In this paper, a strategy is proposed that allows the use of inexact solutions of these subproblems, which, as proved in the paper, can be incorporated without the loss of theoretical convergence guarantees.  Numerical experiments show that by exploiting inexact subproblem solutions, one can consistently reduce the computational effort required by a GS method.  Additionally, a strategy is proposed for aggregating gradient information after a subproblem is solved (potentially inexactly), as has been exploited in bundle methods for nonsmooth optimization.  It is proved that the aggregation scheme can be introduced without the loss of theoretical convergence guarantees.  Numerical experiments show that incorporating this gradient aggregation approach substantially reduces the computational effort required by a GS method.
\end{abstract}


\section{Introduction}

The gradient sampling (GS) methodology has proved to be effective for solving nonsmooth, nonconvex minimization problems.  Based on the conceptually simple idea of computing an approximate $\epsilon$-steepest-descent direction at a point by finding the minimum-norm element of the convex hull of gradients evaluated at randomly generated nearby points, one can prove convergence to stationarity of a GS method under relatively loose assumptions.  That said, here are two ways in which implementations of GS methods could be more efficient:
\begin{hitemize}
  \item Each iteration of a GS method requires the solution of a convex quadratic subproblem (QP) for computing a search direction.  The overall computational expense of a GS method could be reduced if one could terminate each call to a QP solver early, then employ the inexact QP solution as the search direction in the ``outer'' GS method.  Such an inexact solution might cause a search direction to be less productive than if an exact QP solution were computed, meaning that more ``outer'' iterations may be required.  However, as in other optimization algorithms that exploit inexact subproblem solutions, one might still obtain overall computational savings through consistently reduced per-iteration costs.
  \item Bundle methods represent another important class of algorithms for nonsmooth minimization.  Implementations of bundle methods can be made significantly more efficient through the use of \emph{subgradient aggregation}, wherein one can compress the information from a QP solution such that a subsequent QP can be solved more rapidly.  Implementations of GS methods could be made more efficient if such an idea could be incorporated.
\end{hitemize}

In this paper, we propose enhancements to the GS methodology such that one can exploit inexact subproblem solutions and \emph{gradient} aggregation.  (We do not refer to \emph{subgradient} aggregation since the GS methodology requires the identification of points at which the objective function is differentiable or even continuously differentiable, at which gradients are to be evaluated, when search directions are being computed.)  We show techniques for exploiting these ideas that maintain the convergence guarantees of previously proposed GS methods.  Implementations of our ideas in a C++ software package show that exploiting both inexact subproblem solutions and gradient aggregation can lead to consistently noticable reductions in required computational effort.

\subsection{Literature Review}

The GS methodology was introduced by~\cite{BurkLewiOver05}; see also \cite{BurkLewiOver02}.  Shortly after, \cite{Kiwi07} showed an elegant convergence analysis of a GS method, and showed how the convergence guarantees could be maintained by multiple variations of the basic approach.  Later, \cite{CurtQue13,CurtQue15} showed how one could sample gradients adaptively and introduce quasi-Newton Hessian approximations to maintain convergence guarantees while improving practical performance.  (Here, as is common in the literature on quasi-Newton methods for solving nonsmooth optimization problems, we use the term ``Hessian approximation'' loosely; rather than as an approximate second-derivative matrix, it should merely be thought of as a matrix that approximates local changes in the gradient at points at which $f$ is differentiable.)  See also \cite{CurtRobiZhou19} for how to loosen the restrictions on the Hessian approximation scheme.  A feature of the algorithms in all of these articles is that the analyses require that the convex QP subproblems for computing search directions must be solved exactly in every iteration.

A method for reducing the costs associated with solving QPs in a GS method was proposed by~\cite{MortMost19}.  In this work, the authors argue that an ``ideal'' direction, which can be computed using a relatively inexpensive procedure, can be used in place of a QP solution when it is found to be sufficiently large in norm.  The authors argue that convergence guarantees are maintained with this replacement, and show empirically that fewer QPs need to be solved.  Our proposed approach is different from this one in two main respects.  First, rather than prescribe a formula for a particular direction that may be used, our algorithm involves conditions for an inexact QP solution that are more generic.  This gives more computational flexibility to the algorithm.  Second, whereas the algorithm by \cite{MortMost19} still requires that some QPs be solved exactly---such as when the ``ideal'' direction is too small in norm, which occurs when approaching stationarity---our algorithm allows for inexact solutions of the QPs in all cases.

GS ideas have been extended in various ways, such as to attain good local convergence rate properties \citep{HeloSantSimo17} and to solve constrained optimization problems \citep{CurtOver12,HossUsch17,TangLiuJianLi14}.  See \cite{BurkCurtLewiOverSimo19} for further discussion.  Such extensions are beyond the scope of this article, wherein we focus on techniques for unconstrained minimization that ensure convergence from an arbitrarily chosen starting point.  That said, our proposed enhancements could be employed in conjunction with these extensions.

Another prevailing methodology for solving nonsmooth optimization problems is the class of bundle methods, which have a long history \citep{ApkaNollProt08,Kiwi85,Kiwi85b,Kiwi96,HaarMietMaek04,HaarMietMaek07,HareSaga10,HiriLema93b,LemaNemiNest95,LuksVlce98,Miff82,MiffSaga05,Rusz06,SchrZowe92}.  The technique employed in some bundle methods that is relevant for this paper is that of subgradient aggregation; see, e.g., \cite{Kiwi85}.  The use of aggregation in this paper is similar, although the surrounding convergence analysis is different due to the distinct differences in the convergence analyses of bundle and GS methodologies.  For one thing, convergence analyses of GS methods are inherently probabilistic due to the random sampling of points.

\subsection{Notation}

We write $\R{}$ to denote the set of real numbers, $\R{n}$ to denote the set of $n$-dimensional real vectors, and $\R{m \times n}$ to denote the set of $m$-by-$n$-dimensional real matrices.  We write $\N{} := \{0,1,2,\dots\}$ to denote the set of nonnegative integers and use $\mathds{1}$ to denote a vector of ones whose length is determined by the context in which it appears (e.g., through an inner product with a vector of known length).

Throughout the paper, we consider the minimization problem
\bequation\label{prob.f}
  \min_{x \in \R{n}}\ f(x),
\eequation
where the objective function $f : \R{n} \to \R{}$ satisfies the following assumption.
\bassumption\label{ass.f}
  The objective function $f$ is bounded below over $\R{n}$, locally Lipschitz on $\R{n}$, and continuously differentiable on an open set $\Dcal$ with full measure in $\R{n}$.
\eassumption

We propose GS algorithms, each of which is designed to produce an iterate sequence---i.e., $\{x_k\}$ with $x_k \in \R{n}$ for all $k \in \N{}$---converging to stationarity of $f$, which is to say that any cluster point of $\{x_k\}$ is stationary for~$f$.  Throughout, we refer to stationarity in the sense of~\cite{Clar83}.  Such stationarity for $f$ can be defined as follows.  By Rademacher's theorem under Assumption~\ref{ass.f}, the (Clarke) set of generalized gradients of $f$ at $x \in \R{n}$ is given by
\bequation\label{eq.subdiff}
  \eth f(x) = \conv\left\{ \lim_{k\to\infty} \nabla f(x_k) : \text{$\{x_k\} \to x$ with $x_k \in \Dcal$ for all $k \in \N{}$}\right\};
\eequation
see \cite[Theorem~2.5.1]{Clar83}.  For $\epsilon \in [0,\infty)$, the set of $\epsilon$-generalized gradients of $f$ at $x \in \R{n}$ is
\bequation\label{eq.subdiff_eps}
  \eth_\epsilon f(x) = \conv \eth f(\Bmbb(x,\epsilon)),\ \ \text{where}\ \ \Bmbb(x,\epsilon) := \{\xbar \in \R{n} : \|\xbar - x\|_2 \leq \epsilon\}.
\eequation
One finds that $\eth_0 f(x) \equiv \eth f(x)$; see \citep[Corollary~2.5]{Gold77}.  A point $x \in \R{n}$ is said to be $\epsilon$-stationary for $f$ if $0 \in \eth_\epsilon f(x)$ and is said to be stationary for $f$ if $0 \in \eth f(x)$.

The first algorithm that we propose (see Algorithm~\ref{alg.gs_inexact} on page~\pageref{alg.gs_inexact}) has a nested loop (with the ``inner'' loop being stated in Algorithm~\ref{alg.gs_inexact_step} on page~\pageref{alg.gs_inexact_step}).  Iterations for the ``outer'' loop are indexed by $k \in \N{}$.  We apply this iteration number subscript to other values---in addition to $x_k$---computed in the outer loop of the algorithm.  Iterations for the inner loop are indexed by $j \in \N{}$.  Quantities computed during the inner loop are denoted with a double-subscript; e.g., $d_{k,j}$.

One could remove from Assumption~\ref{ass.f} the assumption that~$f$ is bounded below, in which case the methods that we propose would terminate finitely at a stationary point for $f$ or, with probability one, generate iterates that either converge to ($\epsilon$-)stationarity for $f$ (see Theorems~\ref{th.gs_inexact} and \ref{th.gs_aggregation}) or have objective values that diverge to $-\infty$.  However, to focus on the more interesting setting, we include in Assumption~\ref{ass.f} that $f$ is bounded below, meaning this latter case cannot occur.  For the algorithms that we propose to be well posed, one only needs to assume that $f$ is (not necessarily continuously) differentiable in an open set with full measure in $\R{n}$.  However, a theoretical guarantee of convergence to stationarity requires that $f$ be continuously differentiable over such a set, as we have included in Assumption~\ref{ass.f}.  See \cite{BurkCurtLewiOverSimo19} for further discussion.

\subsection{Outline}

In~Section~\ref{sec.inexact}, we propose and analyze an algorithm that employs inexact subproblem solutions.  In Section~\ref{sec.aggregation}, we propose gradient aggregation within a GS method and show that it can be used while maintaining the same guarantees as the method from Section~\ref{sec.inexact}.  Numerical experiments employing both techniques are presented in Section~\ref{sec.numerical}.  Concluding remarks are given in Section~\ref{sec.conclusion}.

\section{GS Algorithm with Inexact Subproblem Solutions}\label{sec.inexact}

We propose a GS algorithm that allows for the use of inexact subproblem solutions in each iteration.  In this section, we present the proposed algorithm, then prove that iterates generated by it converge to ($\epsilon$-)stationarity with probability one.  In our presentation, we focus on the components of the algorithm and analysis that are distinct from previous GS methods.  Components that are not unique are provided in an online companion to this article.

\subsection{Algorithm Description}

In iteration $k \in \N{}$ of our proposed algorithm, an iterate $x_k \in \Dcal$ is available along with a sampling radius $\epsilon_k \in (0,\infty)$, a set of sample points
\bequationNN
  \text{$\Xcal_k := \{x_{k,0},x_{k,1},\dots,x_{k,p_k}\} \subset \Bmbb(x_k,\epsilon_k) \cap \Dcal$ where $x_{k,0} \equiv x_k$ for some $p_k \in \N{}$},
\eequationNN
and the corresponding matrix of gradients
\bequation\label{eq.G}
  G_k := \bbmatrix \nabla f(x_{k,0}) & \nabla f(x_{k,1}) & \cdots & \nabla f(x_{k,p_k}) \ebmatrix \in \R{n \times (p_k+1)}.
\eequation
Given this matrix of gradients, a symmetric positive definite Hessian approximation $H_k$, and the corresponding inverse $W_k := H_k^{-1}$, the search direction is computed by \emph{approximately} solving the primal-dual pair of quadratic optimization problems (QPs) given by
\bequation\label{prob.qp}
  (P) := \left\{
  \baligned
  \min_{(d,z)\in\R{n+1}} &\ z + \thalf \|d\|_{H_k}^2 \\
  \st &\ G_k^Td \leq z \mathds{1}
  \ealigned
  \right\}
  \ \ \text{and}\ \ 
  (D) := \left\{
  \baligned
  \max_{y\in\R{p_k+1}} &\ -\thalf \|G_ky\|_{W_k}^2 \\
  \st &\ \mathds{1}^Ty = 1,\ y \geq 0
  \ealigned
  \right\}.
\eequation
We assume that both $H_k$ and $W_k$ are available for all $k \in \N{}$.  It is straightforward to maintain both approximations through the use of quasi-Newton techniques.

Letting $(d_{k,*},z_{k,*})$ denote the optimal solution of $(P)$ for each $k \in \N{}$, one finds that the solution component $d_{k,*}$ can be viewed as the minimizer of the piecewise quadratic function
\bequationNN
  \max_{i \in \{0,\dots,p_k\}}\ \{\nabla f(x_{k,i})^T d\} + \thalf \|d\|_{H_k}^2.
\eequationNN
The optimal solution $y_{k,*}$ of $(D)$ can be viewed as the vector such that $G_ky_{k,*}$ is the least $W_k$-norm element of the convex hull of the columns of $G_k$, i.e., the $W_k$-projection of the origin onto this hull.  The following lemma reveals important properties of these solutions.

\blemma\label{lem.exact_solution}
  For all $k \in \N{}$, either $(d_{k,*},z_{k,*}) = (0,0)$ and the origin lies in the convex hull of the columns of $G_k$, or $d_{k,*}$ is a direction of strict descent for $f$ at $x_k$ with
  \bequation\label{eq.descent}
  \nabla f(x_k)^Td_{k,*} \leq -d_{k,*}^TH_kd_{k,*} < 0.
  \eequation
  In all cases, $d_{k,*} = -W_kG_ky_{k,*}$ and $\|G_ky_{k,*}\|_{W_k} = \|d_{k,*}\|_{H_k}$.
\elemma
\proof{Proof.}
  The properties follow from the Karush-Kuhn-Tucker conditions for \eqref{prob.qp}; see, e.g., \cite[Eq.~(27)]{CurtQue13} and \cite[Lemma~2.2]{CurtQue15}. \Halmos
\endproof

As our focus is on an algorithm that solves \eqref{prob.qp} \emph{approximately} for all $k \in \N{}$, the statement of our algorithm is facilitated by defining, in each ``outer'' iteration, sequences of inner iterates of a solver for the primal-dual subproblems~\eqref{prob.qp}.  Let $\{(d_{k,j},z_{k,j})\}$ and $\{y_{k,j}\}$ be sequences of primal and dual iterates, respectively, generated when \eqref{prob.qp} is solved iteratively.  Our algorithm requires that \emph{both} primal and dual QP iterate sequences are generated.  However, this is not an expensive requirement.  After all, motivated by Lemma~\ref{lem.exact_solution}, one may choose for a given $y_{k,j} \in \R{p_k+1}$ to set
\bequation\label{eq.primal_from_dual}
  d_{k,j} \gets -W_kG_ky_{k,j}\ \ \text{and}\ \ z_{k,j} \gets \max_{i\in\{0,\dots,p_k\}}\ \nabla f(x_{k,i})^Td_{k,j},
\eequation
in which case one only needs to generate a dual iterate sequence and a corresponding sequence of primal-feasible solutions is obtained through \eqref{eq.primal_from_dual}.  In addition, to reduce expense, one does not need to evaluate \eqref{eq.primal_from_dual} in each inner iteration; one might only evaluate it and check for termination periodically and/or after an initial number of inner iterations.  In any case, for the sake of generality, we define our algorithm to allow $\{d_{k,j}\}$ and $\{-W_kG_ky_{k,j}\}$ to differ.

With respect to the QP solver, we merely assume that the following holds.

\bassumption\label{ass.qp}
  For all $k \in \N{}$, the primal and dual iterates when solving \eqref{prob.qp} satisfy $\{(d_{k,j},z_{k,j},y_{k,j})\} \to (d_{k,*},z_{k,*},y_{k,*})$.  In addition, for all $(k,j) \in \N{} \times \N{}$, one has
  \bequationNN
  G_k^Td_{k,j} \leq z_{k,j} \mathds{1},\ \ \mathds{1}^Ty_{k,j} = 1,\ \ \text{and}\ \ y_{k,j} \geq 0,
  \eequationNN
  i.e., $(d_{k,j},z_{k,j},y_{k,j})$ is primal-dual feasible for all $(k,j) \in \N{} \times \N{}$.
\eassumption

Under Assumption~\ref{ass.qp}, the primal and dual iterates satisfy weak duality with respect to \eqref{prob.qp} for all $(k,j) \in \N{} \times \N{}$.  In particular, defining the QP primal and dual objective functions $q_k : \R{n} \times \R{} \to \R{}$ and $\theta_k : \R{p_k+1} \to \R{}$, respectively, where
\bequationNN
  q_k(d,z) = z + \thalf \|d\|_{H_k}^2\ \ \text{and}\ \ \theta_k(y) = -\thalf \|G_k y\|_{W_k}^2,
\eequationNN
one has that $q_k(d_{k,j},z_{k,j}) \geq \theta_k(y_{k,j})$ for all $(k,j) \in \N{} \times \N{}$.

Our algorithm with inexact subproblem solutions is stated as Algorithm~\ref{alg.gs_inexact} on page~\pageref{alg.gs_inexact}, for which the details of the search direction computation are stated in Algorithm~\ref{alg.gs_inexact_step} on page~\pageref{alg.gs_inexact_step}.  The statement of Algorithm~\ref{alg.gs_inexact} focuses on its unique aspects related to the conditions that we require of inexact QP solutions.  Other subroutines that we employ for the line search, iterate perturbation strategy (a feature required by the theoretical convergence analyses of all GS methods), sample set updates, and quasi-Newton updates are similar to those used in \citep{CurtQue15,CurtRobiZhou19}.  Hence, we relegate them to the online companion.  The algorithm also requires a subroutine for setting parameters related to the quasi-Newton updates that influence the line search subroutine.  The approach is derived from properties of the self-correcting nature of quasi-Newton updating; see \cite{ByrdNoce89}.  This subroutine is also provided in the online companion.

Each iteration of Algorithm~\ref{alg.gs_inexact_step} takes a new approximate subproblem solution from the QP solver.  The loop terminates in one of two situations.  If \eqref{eq.eps_update} holds, then one has obtained a dual iterate such that the corresponding convex combination of columns of~$G_k$ is sufficiently small in appropriate norms.  In this case, one has identified that the current iterate may be sufficiently close to $\epsilon_k$-stationarity, in which case Algorithm~\ref{alg.gs_inexact} reduces the sampling radius.  On the other hand, if \eqref{eq.inexact_step_basic} holds along with either \eqref{eq.inexact_step_1} or \eqref{eq.inexact_step_2}, then our analysis in the following subsection reveals that a sufficiently accurate QP solution yielding a direction of sufficient descent has been obtained.  The condition \eqref{eq.inexact_step_basic} is motivated by Lemma~\ref{lem.exact_solution}, specifically \eqref{eq.descent}, since $(d_k,y_k) = (-W_kG_ky_{k,j_\theta},y_{k,j_\theta})$ yields
\bequationNN
  \nabla f(x_k)^Td_k = -\nabla f(x_k)^TW_kG_ky_k \stackrel{\eqref{eq.inexact_step_basic}}{\leq} -\kappa y_k^TG_k^TW_kG_ky_k = -\kappa d_k^TH_kd_k.
\eequationNN
The role played by conditions \eqref{eq.inexact_step_1} and \eqref{eq.inexact_step_2}, which make use of the values defined in \eqref{eq.tau1} and \eqref{eq.tau2}, is explained in the following subsection.

Notice that an implementation of Algorithm~\ref{alg.gs_inexact_step} does not require storage and a search through all previous subproblem solutions when determining the indices in Line~\ref{step.best}.  One only needs to store the best (in terms of objective values) primal and dual solution estimates during the loop and employ these values when checking for termination of the loop.  Line~\ref{step.best} is only written in this manner for ease of exposition, and to allow us to consider situations in which these inner iterations do not necessarily produce primal and dual subproblem solutions that have objective values that converge monotonically to the optimal value.

\begin{algorithm}[ht]
  \caption{Gradient Sampling Algorithm with Inexact Subproblem Solutions}
  \label{alg.gs_inexact}
  \begin{algorithmic}[1]
  \Require $(\sigma,\underline\alpha) \in (0,\infty)^2$; $(\iota,\underline\phi) \in (0,1)^2$; $\overline\phi \in (1,\infty)$; $\psi \in (0,1]$; $p \in \N{}$ with $p \geq n+1$; $x_0 \in \Dcal$; $H_0 \succ 0$; $\epsilon_0 \in (0,\infty)$.
  \State Set $W_0 \gets H_0^{-1}$, $\Xcal_0 \gets \{x_0\}$, $p_0 \gets 0$, $G_0$ by \eqref{eq.G}, and $\sigma_0 \gets \sigma$.
  \State Set $(\underline\eta,\mu) \in (0,1) \times (1,\infty)$ by Algorithm~\ref{alg.sufficient_decrease}.
  \State Choose $\overline\eta \in (\underline\eta,1)$.
  \For{\textbf{all} $k \in \N{}$}
  \If{$\|\nabla f(x_k)\|_2=0$}\label{step.gs_inexact.terminate}
  \State \textbf{terminate} and \textbf{return} the stationary point $x_k$.
  \EndIf
  \State Set $y_{k,j_\theta}$ by Algorithm~\ref{alg.gs_inexact_step} (page~\pageref{alg.gs_inexact_step}).
  \State Set $(d_k,y_k) \gets (-W_kG_ky_{k,j_\theta},y_{k,j_\theta})$.
  \State Set $\alpha_k \geq 0$ by Algorithm~\ref{alg.ls}. \label{step.ls}
  \State \textbf{if} \eqref{eq.eps_update} holds (with $y_{k,j_\theta} \equiv y_k$) \label{step.eps_decrease}
  \State \hskip1.5em set $\epsilon_{k+1} \gets \psi \epsilon_k$ and $\sigma_{k+1} \gets \sigma$;
  \State \textbf{else if} $\alpha_k \geq \underline\alpha$
  \State \hskip1.5em set $\epsilon_{k+1} \gets \epsilon_k$ and $\sigma_{k+1} \gets \sigma_k$;
  \State \textbf{else}
  \State \hskip1.5em set $\epsilon_{k+1} \gets \epsilon_k$ and $\sigma_{k+1} \gets \iota \sigma_k$.
  \State \textbf{end if}
  \State Set $x_{k+1} \in \Dcal$ by Algorithm~\ref{alg.perturb}. \label{step.update}
  \State Set $(H_{k+1},W_{k+1})$ by Algorithm~\ref{alg.hessian}.\label{step.HW}
  \State Set $(\Xcal_{k+1},p_{k+1})$ by Algorithm~\ref{alg.sample} and $G_{k+1}$ by \eqref{eq.G}. \label{step.sample}
  \EndFor
  \end{algorithmic}
\end{algorithm}

\begin{algorithm}[ht]
  \caption{Search Direction Computation for Algorithm~\ref{alg.gs_inexact}}
  \label{alg.gs_inexact_step}
  \begin{algorithmic}[1]
  \Require $\nu \in (0,\infty)$; $(\rho,\kappa) \in (0,1)^2$
  \State Set
  \bequation\label{eq.tau1}
  \tau_k \gets \sigma_k^2 + 2\sigma_k \in (0,1).
  \eequation
  \For{\textbf{all} $j \in \N{}$}
  \State Set $j_q \gets \displaystyle \argmin_{i \in \{0,\dots,j\}}\ q_k(d_{k,i},z_{k,i})$ and $j_\theta \gets \displaystyle \argmax_{i \in \{0,\dots,j\}}\ \theta_k(y_{k,i})$. \label{step.best}
  \State \textbf{if} $q_k(d_{k,j_q},z_{k,j_q}) \geq 0$, \textbf{then} set $\lambda_{k,j_q} = \infty$
  \State \textbf{else} set
  \bequation\label{eq.tau2}
  \lambda_{k,j_q} \gets \max\left\{1-\frac{\sigma_k^2+2\sigma_k}{\frac{\theta_k(y_{k,0})}{q_k(d_{k,j_q},z_{k,j_q})}-1},\rho\right\}.
  \eequation
  \State \textbf{end if}
  \State \textbf{if}
  \bequation\label{eq.eps_update}
  \max\{\|W_kG_ky_{k,j_\theta}\|_2,\|G_ky_{k,j_\theta}\|_2\} \leq \nu \epsilon_k,
  \eequation
  \State \hskip1.5em \textbf{then} \textbf{terminate} and \textbf{return} $y_{k,j_\theta}$; \label{step.bfgsgs_qp1}
  \State \textbf{else if}
  \bequation\label{eq.inexact_step_basic}
  -\nabla f(x_k)^TW_kG_ky_{k,j_\theta} \leq -\kappa y_{k,j_\theta}^TG_k^TW_kG_ky_{k,j_\theta}
  \eequation
  \State \hskip1.5em \textbf{and either}
  \bequation\label{eq.inexact_step_1}
  q_k(d_{k,j_q},z_{k,j_q}) - \theta_k(y_{k,j_\theta}) \leq \tau_k (- q_k(d_{k,j_q},z_{k,j_q}))
  \eequation
  \State \hskip1.5em \textbf{or}
  \bequation\label{eq.inexact_step_2}
  \theta_k(y_{k,j_\theta}) - \theta_k(y_{k,0}) \geq \lambda_{k,j_q} (q_k(d_{k,j_q},z_{k,j_q}) - \theta_k(y_{k,0}))
  \eequation
  \State \hskip1.5em \textbf{then} \textbf{terminate} and \textbf{return} $y_{k,j_\theta}$.
  \State \textbf{end if}
  \EndFor
  \end{algorithmic}
\end{algorithm}

\subsection{Inexactness conditions for the QP solver}

Convergence analyses of GS methods rely on a fundamental property of any compact, convex set, call it $\Scal \subseteq \R{n}$, that does not contain the origin.  Intuitively, this property is that if $u \in \Scal$ is sufficiently close to the projection of the origin onto $\Scal$, then $u$ makes a sufficiently acute angle (with respect to a given metric) with any $v \in \Scal$.  Such a lemma appears as \cite[Lemma~3.1]{BurkLewiOver05} and \cite[Lemma~3.1]{Kiwi07}, and is proved in a more general setting as \cite[Lemma~3.5]{CurtQue15}.  The conditions that we impose on inexact subproblem solutions are motivated by trying to ensure a property of this type, but in an even slightly more general setting.  Specifically, the lemma that we use is the following.  In the lemma, we refer to the concept of a $W$-projection (with $W \succ 0$) of the origin onto a compact, convex set $\Scal$, i.e.,
\bequation\label{eq.projection}
  P_W(\Scal) := \argmin_{s \in \Scal} \|s\|_W.
\eequation
Our new generalization of the lemma can be seen in the inequality \eqref{eq.angle_lemma}, which does not require that a given vector $u \in \Scal$ is sufficiently close to the $W$-projection of the origin, but merely sufficiently close to a small enough neighborhood of this projection.

\blemma\label{lem.angle_lemma}
  Suppose $\Scal \subseteq \R{n}$ is a compact and convex set with $0 \notin \Scal$.  For any $\beta \in(0,1)$ and $W \succ 0$, there exists $(\varsigma,\delta) \in (0,\infty)^2$ such that, for any $(u,v,\bar\varsigma,\bar\delta) \in \Scal \times \Scal \times (0,\varsigma] \times (0,\delta]$ with
  \bequation\label{eq.angle_lemma}
  \|u\|_W \leq (1 + \varsigma)\|P_W(S)\|_W + \delta
  \eequation
  $($where $P_W(\Scal)$ is defined in \eqref{eq.projection}$)$, it follows that $v^TWu > \beta \|u\|_W^2$.
\elemma
\proof{Proof.}
  Consider arbitrary $\beta \in (0,1)$ and $W \succ 0$.  To derive a contradiction, suppose the implication is false, that is, for any $(\varsigma,\delta) \in (0,\infty)^2$ there exists $(u,v,\bar\varsigma,\bar\delta) \in \Scal \times \Scal \times (0,\varsigma] \times (0,\delta]$ with
  \bequationNN
    \|u\|_W \leq (1+\bar\varsigma)\|P_W(\Scal)\|_W + \bar\delta\ \ \text{and}\ \ v^TWu \leq \beta \|u\|_W^2.
  \eequationNN
  This means that there exist infinite sequences $\{u_i\} \subset \Scal$ and $\{v_i\} \subset \Scal$ such that
  \bequation\label{eq.u_i}
    \|u_i\|_W \leq (1 + 1/i)\|P_W(\Scal)\|_W + 1/i\ \ \text{and}\ \ v_i^TWu_i \leq \beta \|u_i\|_W^2\ \ \text{for all}\ \ i \in \N{}.
  \eequation
  Since $\Scal$ is compact, these sequences have convergent subsequences; hence, without loss of generality, one can assume that $\{u_i\} \to u$ and $\{v_i\} \to v$ for some $(u,v) \in \Scal \times \Scal$ with
  \bequation\label{eq.jobin}
    v^TWu \leq \beta \|u\|_W^2.
  \eequation
  On the other hand, by the definition of $\{u_i\}$, it follows that $u = P_W(\Scal)$, which is nonzero since $\Scal$ does not include the origin.  Moreover, by \cite[Proposition~1.1.8]{Bert09} and the definition of $P_W(\Scal)$ (as the $W$-projection of the origin onto $\Scal$), one finds that
  \bequationNN
    (0 - u)^TW(v - u) \leq 0 \iff v^TWu \geq \|u\|_W^2,
  \eequationNN
  which contradicts \eqref{eq.jobin} since $\beta \in (0,1)$. \Halmos
\endproof

Our goal now is to prove two lemmas that motivate the use of \eqref{eq.inexact_step_1} and \eqref{eq.inexact_step_2} as stopping conditions for the loop in Algorithm~\ref{alg.gs_inexact_step}.  Specifically, if $\theta_k(y_{k,*}) < 0$, each lemma shows that these conditions---\eqref{eq.inexact_step_1} and \eqref{eq.inexact_step_2}, respectively, in the two lemmas---imply that
\bequation\label{eq.theta_good}
  \baligned
  0 &> \theta_k(y_{k,j_\theta}) \geq (1 + \sigma_k)^2\theta_k(y_{k,*}) \\
  \iff 0 &< \|G_ky_{k,j_\theta}\|_{W_k} \leq (1 + \sigma_k) \|G_ky_{k,*}\|_{W_k}.
  \ealigned
\eequation
Importantly, these algorithmic conditions imply that \eqref{eq.theta_good} holds \emph{without knowledge of~$y_{k,*}$}.  The inequalities in \eqref{eq.theta_good} are important since they, along with Lemma~\ref{lem.angle_lemma} (c.f.~\eqref{eq.angle_lemma}), play a central role in our convergence analysis in \S\ref{sec.convergence} for Algorithm~\ref{alg.gs_inexact}.

\blemma\label{lem.approx}
  Suppose that, in iteration $k \in \N{}$ of Algorithm~\ref{alg.gs_inexact}, one has $\theta_k(y_{k,*}) < 0$.  In addition, suppose that, during iteration $j \in \N{}$ of Algorithm~\ref{alg.gs_inexact_step} $($during outer iteration $k \in \N{}$$)$, one finds that \eqref{eq.inexact_step_1} holds.  Then, \eqref{eq.theta_good} holds.
\elemma
\proof{Proof.}
  By weak duality for \eqref{prob.qp}, one has that
  \bequationNN
  \baligned
  \theta_k(y_{k,*}) - \theta_k(y_{k,j_\theta}) &\leq q_k(d_{k,j_q},z_{k,j_q}) - \theta_k(y_{k,j_\theta}) \\ \text{and}\ \ -q_k(d_{k,j_q},z_{k,j_q}) &\leq -\theta_k(y_{k,*}).
  \ealigned
  \eequationNN
  Combined with \eqref{eq.inexact_step_1} and \eqref{eq.tau1}, it follows that
  \bequationNN
  \theta_k(y_{k,*}) - \theta_k(y_{k,j_\theta}) \leq \tau_k(-\theta_k(y_{k,*})) = (\sigma_k^2 + 2\sigma_k)(-\theta_k(y_{k,*})),
  \eequationNN
  which shows that \eqref{eq.theta_good} holds, as desired. \Halmos
\endproof

When $\theta_k(y_{k,*}) < 0$, weak duality for \eqref{prob.qp} implies that \eqref{eq.inexact_step_1} can hold only if $q_k(d_{k,j_q},z_{k,j_q}) < 0$.  Hence, one does not need to check if $q_k(d_{k,j_q},z_{k,j_q}) < 0$ holds before employing \eqref{eq.inexact_step_1} as a stopping condition for the QP solver.  By contrast, the next lemma shows that \eqref{eq.inexact_step_2} should be used as a stopping condition for the QP solver only if $q_k(d_{k,j_q},z_{k,j_q}) < 0$.  Algorithm~\ref{alg.gs_inexact_step} ensures this by setting $\lambda_{k,j_q} \gets \infty$ when $q_k(d_{k,j_q},z_{k,j_q}) \geq 0$, and otherwise the lemma shows that $\lambda_{k,j_q} \in (0,1)$.

\blemma\label{lem.approx_2}
  Suppose that, in iteration $k \in \N{}$ of Algorithm~\ref{alg.gs_inexact}, one has $\theta_k(y_{k,*}) < 0$.  In addition, suppose that, during iteration $j \in \N{}$ of Algorithm~\ref{alg.gs_inexact_step} $($during outer iteration $k \in \N{}$$)$, one finds that $q_k(d_{k,j_q},z_{k,j_q}) < 0$ and \eqref{eq.inexact_step_2} holds.  Then, \eqref{eq.theta_good} holds.
\elemma
\proof{Proof.}
  By $q_k(d_{k,j_q},z_{k,j_q}) < 0$ and weak duality for \eqref{prob.qp}, one finds in \eqref{eq.tau2} that
  \bequation\label{eq.lambda_bound}
  \frac{\theta_k(y_{k,0})}{q_k(d_{k,j_q},z_{k,j_q})} \geq 1.
  \eequation
  If $\theta_k(y_{k,0}) = q_k(d_{k,j_q},z_{k,j_q})$, then $(d_{k,j_q},z_{k,j_q},y_{k,0})$ is a primal-dual solution of \eqref{prob.qp} and $\theta_k(y_{k,0}) = \theta_k(y_{k,j_\theta}) = \theta_k(y_{k,*})$, which means that \eqref{eq.theta_good} holds.  Hence, we may proceed under the assumption that $\theta_k(y_{k,0}) < q_k(d_{k,j_q},z_{k,j_q}) < 0$, which implies that \eqref{eq.lambda_bound} holds strictly.  Observing \eqref{eq.tau2}, one finds $\lambda_{k,j_q} \in (0,1)$.  This fact, \eqref{eq.inexact_step_2}, and weak duality for \eqref{prob.qp} imply
  \bequationNN
  \baligned
  \theta_k(y_{k,j_\theta}) - \theta_k(y_{k,0})
  &\geq \lambda_{k,j_q}(q_k(d_{k,j_q},z_{k,j_q}) - \theta(y_{k,0})) \\
  &\geq \lambda_{k,j_q}(\theta_k(y_{k,*}) - \theta(y_{k,0})) \geq 0,
  \ealigned
  \eequationNN
  which along with $\lambda_{k,j_q} \in (0,1)$ and the facts that $\theta_k(y_{k,*}) < 0$ and
  \bequationNN
  \frac{\theta_k(y_{k,0})}{\theta_k(y_{k,*})} \leq \frac{\theta_k(y_{k,0})}{q_k(d_{k,j_q},z_{k,j_q})} \implies \theta_k(y_{k,0}) \geq \frac{\theta_k(y_{k,0})}{q_k(d_{k,j_q},z_{k,j_q})}\theta_k(y_{k,*})
  \eequationNN
  implies that
  \bequation\label{eq.dual_lower}
  \baligned
  \theta_k(y_{k,j_\theta})
  &\geq \lambda_{k,j_q} \theta_k(y_{k,*}) + (1 - \lambda_{k,j_q})\theta_k(y_{k,0}) \\
  &\geq \(\lambda_{k,j_q} + (1 - \lambda_{k,j_q})\frac{\theta_k(y_{k,0})}{q_k(d_{k,j_q},z_{k,j_q})}\)\theta_k(y_{k,*}).
  \ealigned
  \eequation
  In addition, one finds that $\lambda_{k,j_q}$ in \eqref{eq.tau2} satisfies
  \bequationNN
  \lambda_{k,j_q} \geq 1 - \frac{\sigma_k^2+2\sigma_k}{\frac{\theta_k(y_{k,0})}{q_k(d_{k,j_q},z_{k,j_q})}-1} = \frac{\frac{\theta_k(y_{k,0})}{q_k(d_{k,j_q},z_{k,j_q})} - (1 + \sigma_k)^2}{\frac{\theta_k(y_{k,0})}{q_k(d_{k,j_q},z_{k,j_q})} - 1},
  \eequationNN
  implying that
  \bequationNN
  \lambda_{k,j_q} + (1-\lambda_{k,j_q})\frac{\theta_k(y_{k,0})}{q_k(d_{k,j_q},z_{k,j_q})} \leq (1 + \sigma_k)^2,
  \eequationNN
  which along with \eqref{eq.dual_lower} and the fact that $\theta_k(y_{k,*}) < 0$ shows that
  \bequationNN
  \theta_k(y_{k,j_\theta}) \geq \(\lambda_{k,j_q} + (1 - \lambda_{k,j_q})\frac{\theta_k(y_{k,0})}{q_k(d_{k,j_q},z_{k,j_q})}\) \theta_k(y_{k,*}) \geq (1 + \sigma_k)^2 \theta(y_{k,*}),
  \eequationNN
  as desired. \Halmos
\endproof

\subsection{Convergence Analysis}\label{sec.convergence}

In this section, we show under Assumptions~\ref{ass.f} and \ref{ass.qp} that Algorithm~\ref{alg.gs_inexact} either terminates finitely with a stationary point for $f$ or, with probability one, generates a sequence of iterates that converge to stationarity for $f$.  Throughout this section, let $\Kcal$ be the indices of the outer iterations performed by the algorithm before termination (if the algorithm ever terminates) or the failure of a subroutine (if a subroutine ever fails).  The subroutines that may fail are the iteration perturbation procedure (Algorithm~\ref{alg.perturb} in the online companion) and the sample set update (Algorithm~\ref{alg.sample} in the online companion), wherein failure means that a loop does not terminate.  If such an event occurs in iteration $k$, then $\Kcal = \{1,\dots,k\}$.  If the algorithm never terminates and no subroutine ever fails, then one simply has that the iterations performed are $\Kcal = \N{}$.

We begin by showing that the algorithm is well posed along with important properties of the subroutines stated in the online companion.

\blemma\label{lem.well_posed}
  Algorithm~\ref{alg.gs_inexact} is well posed; it either terminates finitely or, with probability one, it performs an infinite number of iterations.  In any case, for all $k \in \Kcal$, the following hold true.
  \benumerate
  \item[(a)] $H_k \succ 0$ and $W_k = H_k^{-1} \succ 0$.
  \item[(b)] $(d_k,y_k)$ satisfies $\|d_k\|_{H_k} = \|G_ky_k\|_{W_k}$.
  \item[(c)] In Line~\ref{step.ls}, Algorithm~\ref{alg.ls} terminates finitely with $\alpha_k \geq 0$.  If $p_k < p$, then $\alpha_k = 0$ or $\alpha_k \in [\underline\alpha,\overline\alpha]$.  Otherwise, if $p_k = p$, then $\alpha_k \in (0,\overline\alpha]$.  In any case, if $\alpha_k > 0$, then 
  \bsubequations\label{eq.wolfe}
    \begin{align}
      &f(x_k) - f(x_k + \alpha_kd_k) > \underline\eta \alpha_k \max\{\|d_k\|_2^2,\|G_ky_k\|_2^2\} \label{eq.armijo} \\
      \text{and}\ \ &v^Td_k \geq \overline\eta \nabla f(x_k)^Td_k,\ \text{where}\ v \in \eth f(x_k+\alpha_kd_k), \label{eq.curv}
    \end{align}
  \esubequations
  or at least \eqref{eq.armijo} holds $($which is sufficient if deemed by Algorithm~\ref{alg.ls}$)$.
  \item[(d)] In Line~\ref{step.update}, Algorithm~\ref{alg.perturb} yields, with probability one, $x_{k+1} \in \Dcal$ satisfying
  \bsubequations\label{eq.wolfe2}
    \begin{align}
      f(x_k) - f(x_{k+1}) &\geq \underline\eta \alpha_k \max\{\|d_k\|_2^2,\|G_ky_k\|_2^2\}, \label{eq.armijo2} \\
      \nabla f(x_{k+1})^Td_k &\geq \overline\eta \nabla f(x_k)^Td_k, \label{eq.curv2} \\
      \text{and}\ \ \|x_k + \alpha_kd_k - x_{k+1}\|_2 &\leq \min\{\alpha_k,\epsilon_k\}\min\{\|d_k\|_2,\|G_ky_k\|_2\}, \label{eq.perturb}
    \end{align}
  \esubequations
  or at least satisfying \eqref{eq.armijo2} and \eqref{eq.perturb} $($which is sufficient if deemed by Algorithm~\ref{alg.perturb}$)$.
  \item[(e)] If Line~\ref{step.sample} is reached and one finds that
  \bequation\label{eq.Gy_more_d_and_alpha_big}
    \|d_k\|_{H_k}^2 \geq \xi \|d_k\|_2^2\ \ \text{and}\ \ \alpha_k \geq \underline\alpha,
  \eequation
  then Algorithm~\ref{alg.sample} yields $\Xcal_{k+1} \gets \{x_{k+1}\}$ and $p_{k+1} \gets 0$; otherwise, with probability one,
  \bequationNN
  \Xcal_{k+1} \gets (\{x_{k+1}\} \cup \Scal_{k+1} \cup (\Xcal_k \cap \Bmbb(x_{k+1},\epsilon_{k+1}))) \subset \Bmbb(x_{k+1},\epsilon_{k+1})\ \text{with}\ p_{k+1} \geq \min\{p_k+1,p\}.
  \eequationNN
  \eenumerate
  Finally, let $\Kcal_{H,W} := \{k \in \Kcal : \alpha_k d_k =: s_k \neq 0\}$, which are the indices for which Algorithm~\ref{alg.hessian} may yield $(H_{k+1},W_{k+1}) \neq (H_k,W_k)$.  If $\Kcal_{H,W}$ is infinite, then for any $\chi \in (0,1)$ there exists $(\underline\mu,\overline\mu) \in (0,\infty)^2$ with $\underline\mu \leq \overline\mu$ such that, for every $K \in \N{}$, the following hold for at least $\lceil \chi K \rceil$ values of $k \in \Kcal_{H,W}$:
  \bsubequations\label{eq.HW}
  \begin{align}
  \underline\mu \|G_ky_k\|_2^2 &\leq \|G_ky_k\|_{W_k}^2 \label{eq.HW_W1} \\
  \text{and}\ \ \|W_kG_ky_k\|_2^2 &\leq \overline\mu \|G_ky_k\|_2^2. \label{eq.HW_W2}
  \end{align}
  \esubequations
  If $\Kcal_{H,W}$ is finite, then such constants exist satisfying \eqref{eq.HW} for all $k \in \Kcal$.
\elemma
\proof{Proof.}
  If the algorithm reaches iteration $k \in \Kcal$ in which the condition in Line~\ref{step.gs_inexact.terminate} holds, then the algorithm terminates finitely.  In this case, all subroutines in iterations $\{0,1,\dots,k-1\}$ must have terminated successfully prior to termination.  Moreover, in this case, \eqref{eq.HW} follows from the fact that only a finite number of iterations are performed and the following proof of part (a) of the lemma:
  \begin{henumerate}
  \item[(a)] The facts that $H_0 \succ 0$ and $W_0 \succ 0$ follow from the initialization of the algorithm.  Now suppose that iteration $1$ is reached.  If $s_0 = 0$, then $H_1 \gets H_0 \succ 0$ and $W_1 \gets W_0 \succ 0$; otherwise, positive definiteness of $H_1$ and $W_1$ follows the fact that \eqref{eq.phi} implies $s_0^Tv_0 > 0$ and from well-known properties of BFGS updating; see, e.g., \cite[Chapter~6]{NoceWrig06}.  Inductively, positive definiteness of $H_k$ and $W_k$ for any $k \in \N{}$ follows by the same arguments.
  \end{henumerate}
  This completes the proof of the lemma for the case when the algorithm reaches $k \in \Kcal$ at which the condition in Line~\ref{step.gs_inexact.terminate} holds.  Hence, we may proceed under the assumption that this condition does not hold for any $k \in \Kcal$.
  
  Suppose that the algorithm reaches iteration $k \in \Kcal$.  To prove that, with probability one, it reaches iteration $k+1$ (i.e., without failure of a subroutine), it suffices to prove parts (b)--(e) (since part (a) has been proved above).
  
  \begin{henumerate}
  \item[(b)] By part (a), one has $H_k \succ 0$ and $W_k \succ 0$, from which it follows that strong duality holds at the primal-dual optimal solution of \eqref{prob.qp}.  Since $\theta_k(y) \leq 0$ for all $y \in \N{}$, there are two cases to consider, namely, whether $\theta_k(y_{k,*}) = 0$ or $\theta_k(y_{k,*}) < 0$.  First, suppose that $\theta_k(y_{k,*}) = 0$.  Since $W_k \succ 0$, this implies that $G_ky_{k,*} = 0$.  Under Assumption~\ref{ass.qp}, we have that $y_{k,j} \to y_{k,*}$.  This limit, the fact that $G_ky_{k,*} = 0$, and the facts that $W_k \succ 0$ and $\epsilon_k > 0$ together imply that \eqref{eq.eps_update} holds for some sufficiently large $j \in \N{}$.  Now suppose that $\theta_k(y_{k,*}) < 0$.  If \eqref{eq.eps_update} holds for any $j \in \N{}$, then the inner loop terminates and there is nothing left to prove; hence, we may proceed assuming that \eqref{eq.eps_update} does not hold for any $j \in \N{}$.  Under Assumption~\ref{ass.qp}, we have that $(d_{k,j},y_{k,j}) \to (d_{k,*},y_{k,*})$.  This limit, continuity of $q_k$ and $\theta_k$, the fact that $\theta_k(y_{k,*}) < 0$, strong duality for \eqref{prob.qp}, Lemma~\ref{lem.exact_solution}, and the fact that $\tau_k \in (0,1)$ imply that \eqref{eq.inexact_step_basic} and \eqref{eq.inexact_step_1} will be satisfied for some sufficiently large $j \in \N{}$.  Finally, the fact that $H_k = W_k^{-1}$ and at termination of the inner loop the algorithm yields $d_k = -W_kG_ky_k$ implies that $\|d_k\|_{H_k} = \|G_ky_k\|_{W_k}$, as desired.
  \item[(c--d)] The proof follows in the same manner as that for \cite[Lemma~2.3]{CurtQue15}.
  \item[(e)] The proof follows in the same manner as that for \cite[Lemma~2.5]{CurtQue15}.
  \end{henumerate}
  Since we have shown that if the algorithm reaches iteration $k \in \Kcal$, then it reaches iteration $k+1$ with probability one, it follows that, again with probability one, an infinite number of iterations are performed.  Finally, with respect to the stated property of the sequence $\{(H_k,W_k)\}_{k \in \Kcal_{H,W}}$, the proof follows in the same manner as that for \cite[Corollary~3.2]{CurtRobiZhou19}. \Halmos
\endproof

The next three lemmas are similar to results previously proved for GS methods.  First, the following lemma is a simple consequence of the previous lemma (specifically, parts~(c) and (e)) and the sample set update strategy, namely, Algorithm~\ref{alg.sample} (in the online companion).  A similar result was proved as \cite[Lemma~3.3]{CurtQue15}.

\blemma\label{lem.alpha}
  If $\Kcal = \N{}$, then $\Kcal_\alpha := \{k \in \N{} : \alpha_k > 0\}$ is infinite.
\elemma
\proof{Proof.}
  Suppose $\Kcal = \N{}$ and observe by Lemma~\ref{lem.well_posed}(c) that $\alpha_k \geq 0$ for all $k \in \N{}$.  In order to derive a contradiction, suppose that there exists an index $k_\alpha \in \N{}$ such that $\alpha_k = 0$ for all $k \in \N{}$ with $k \geq k_\alpha$.  By Lemma~\ref{lem.well_posed}(c), this means that $p_k \leq p-1$ for all $k \geq k_\alpha$.  However, with $\alpha_k = 0$, one finds that \eqref{eq.Gy_more_d_and_alpha_big} does not hold, which by Lemma~\ref{lem.well_posed}(e) implies that $p_{k+1} \geq \min\{p_k+1,p\}$.  This implies the existence of some $k \geq k_\alpha$ such that $p_k \geq p$, which by Lemma~\ref{lem.well_posed}(c) implies that $\alpha_k > 0$, a contradiction of the definition of the index $k_\alpha$. \Halmos
\endproof

The next lemma shows a useful upper bound on the objective function value at iteration $k + 1 \in \Kcal$; for a similar result, see, e.g., \cite[Lemma~3.4]{CurtQue15}.

\blemma\label{lem.inequality}
  If $k+1 \in \Kcal$, then
  \bequationNN
  f(x_{k+1}) \leq f(x_k) - \thalf \underline\eta \|x_{k+1}-x_k\|_2 \max\{\|d_k\|_2,\|G_ky_k\|_2\}.
  \eequationNN
\elemma
\proof{Proof.}
  Suppose $k+1 \in \Kcal$, which implies that $k \in \Kcal$.  Lemma~\ref{lem.well_posed}(a) and (b) imply that $d_k = 0$ if and only if $G_ky_k = 0$.  If $d_k = 0$ and $G_ky_k = 0$, then $x_{k+1} = x_k$ and the result follows trivially.  Otherwise, in iteration $k \in \Kcal$, Lemma~\ref{lem.well_posed}(d) shows that $x_{k+1}$ satisfies \eqref{eq.armijo2} and \eqref{eq.perturb}.  The triangle inequality and~\eqref{eq.perturb} imply
  \bequationNN
  \baligned
  \|x_{k+1} - x_k\|_2
  &\leq \min\{\alpha_k,\epsilon_k\} \min\{\|d_k\|_2,\|G_ky_k\|_2\} + \alpha_k \|d_k\|_2 \\
  &\leq \alpha_k \|d_k\|_2 \min\{2,1 + \|G_ky_k\|_2/\|d_k\|_2\}.
  \ealigned
  \eequationNN
  Hence, along with \eqref{eq.armijo2}, one finds that
  \bequationNN
  \baligned
  f(x_{k+1}) - f(x_k)
  &\leq -\underline\eta \alpha_k \max\{\|d_k\|_2^2,\|G_ky_k\|_2^2\} \\
  &=  -\underline\eta \alpha_k \|d_k\|_2 \max\{\|d_k\|_2,\|G_ky_k\|_2^2/\|d_k\|_2\} \\
  &\leq -\underline\eta \|x_{k+1} - x_k\|_2 \(\frac{\max\{\|d_k\|_2,\|G_ky_k\|_2^2/\|d_k\|_2\}}{\min\{2,1 + \|G_ky_k\|_2/\|d_k\|_2\}}\) \\
  &\leq -\thalf \underline\eta \|x_{k+1} - x_k\|_2 \max\{\|d_k\|_2,\|G_ky_k\|_2\},
  \ealigned
  \eequationNN
  as desired. \Halmos
\endproof

Now we enter the core theory of GS methods.  At its heart is the closure of the convex hull of gradients at points of differentiability in an $\epsilon_k$-neighborhood about a given point $\xbar \in \R{n}$, namely,
\bequation\label{def.Gcal}
  \Gcal(\xbar,\epsilon_k) := \cl \conv \nabla f(\Bmbb(\xbar,\epsilon_k) \cap \Dcal).
\eequation
along with, for any $\omega \in (0,\infty)$, the subset of the Cartesian product of $\epsilon_k$-balls about $x_k$ given by
\bequationNN
  \baligned
  \Tcal_k(\xbar,\omega) := \Bigg\{ & \Xcal_k \in \prod_0^{p_k} (\Bmbb(x_k,\epsilon_k) \cap \Dcal) : \\
  & \|P_{W_k}(\conv(\{\nabla f(x)\}_{x\in \Xcal_k}))\|_{W_k} \leq \|P_{W_k}(\Gcal(\xbar,\epsilon_k))\|_{W_k} + \omega \Bigg\},
  \ealigned
\eequationNN
both of which are defined with respect to each iteration number $k \in \N{}$ and a point $\xbar \in \R{n}$.  (In the definition of $\Tcal_k(\xbar,\omega)$, recall that $P_{W_k}(\cdot)$ has been defined in \eqref{eq.projection}.)  The following lemma, which follows \cite[Lemma 3.2(i)]{Kiwi07}, \cite[Lemma~4.7]{CurtQue13}, and \cite[Lemma~3.6]{CurtQue15}, shows that if the sample set size indictor $p_k$ is sufficiently large and $x_k$ is sufficiently close to $\xbar$, then for any $\omega \in (0,\infty)$ there exists a nonempty open subset of $\Tcal_k(\xbar,\omega)$.  This will be critical in our main result, where we need to show in certain situations that an element of this subset can be found through random sampling of points.

\blemma\label{lem.T}
  Let $\xbar \in \R{n}$ and $\omega \in (0,\infty)$ be given.  If $k \in \Kcal$ and $p_k \geq n+1$, then there exists $\zeta > 0$ such that with $x_k \in \Bmbb(\xbar,\zeta)$ there is a nonempty open $\Tcal \subseteq \Tcal_k(\xbar,\omega)$.
\elemma
\proof{Proof.}
  Using the metric defined by $W_k$, the proof follows the same argument of \cite[Lemma 3.2(i)]{Kiwi07}, which makes use of Carath\'eodory's theorem. \Halmos
\endproof

We now present a convergence theorem for Algorithm~\ref{alg.gs_inexact}.  Much of the proof follows similar arguments as that for \cite[Theorem~3.1]{CurtQue15}, which we present for completeness.  The new features are two-fold: (1) Our algorithm is even less conservative about the Hessian and inverse Hessian updates than the method in \cite{CurtQue15}, so our convergence result relies on arguments about self-correcting properties of BFGS updating that we have stated in Lemma~\ref{lem.well_posed}, which borrows from \cite{CurtRobiZhou19}; and (2) our inexactness conditions and our Lemma~\ref{lem.angle_lemma}, which have not appeared before for GS methods, play critical roles in the proof of the theorem.

\btheorem\label{th.gs_inexact}
  Suppose $\psi \in (0,1)$.  Algorithm~\ref{alg.gs_inexact} either terminates finitely with a stationary point for $f$ or, with probability one, it performs an infinite number of outer iterations.  In the latter case, with probability one, the sampling radius sequence satisfies $\{\epsilon_k\} \searrow 0$ and every cluster point of the iterate sequence $\{x_k\}$ is stationary for~$f$.
\etheorem
\proof{Proof.}
  If Algorithm~\ref{alg.gs_inexact} terminates finitely with a stationary point for $f$, then there is nothing left to prove.  Otherwise, by Lemma~\ref{lem.well_posed}, it follows with probability one that an infinite number of outer iterations are performed, meaning $\Kcal = \N{}$.  Since our desired conclusion only needs to hold with probability one, we may assume going forward that $\Kcal = \N{}$.  Under this assumption, our next aim is to prove that $\{\epsilon_k\} \searrow 0$ with probability one.  We consider two cases.
  
  \begin{hitemize}
  \item \textbf{Case 1:} Suppose that $\Kcal_d := \{k \in \N{} : d_k = 0\}$ is infinite.  By Lemma~\ref{lem.well_posed}(a) and~(b), it follows that $G_ky_k = 0$ for all $k \in \Kcal_d$.  This fact, the fact that $|\Kcal_d| = \infty$, and \eqref{eq.eps_update} imply that $\{\epsilon_k\} \searrow 0$.
  \item \textbf{Case 2:} Suppose that $\Kcal_d := \{k \in \N{} : d_k = 0\}$ is finite.  Let us proceed by supposing that there exists $k_\epsilon \in \N{}$ and a sampling radius $\epsilon \in (0,\infty)$ such that $\epsilon_k = \epsilon$ for all $k \in \N{}$ with $k \geq k_\epsilon$.  Our aim is to show that the existence of such a pair $(k_\epsilon,\epsilon)$ occurs with probability zero.  From \eqref{eq.eps_update},
  \bequation\label{eq.eps_no_update}
  \max\{\|d_k\|_2,\|G_ky_k\|_2\} > \nu \epsilon\ \ \text{for all}\ \ k \geq k_\epsilon.
  \eequation
  On the other hand, Assumption~\ref{ass.f}, Lemma~\ref{lem.inequality}, and \eqref{eq.armijo2} imply that
  \bsubequations\label{eq.sums}
  \begin{align}
  \sum_{k = k_\epsilon}^\infty \|x_{k+1} - x_k\|_2 \max\{\|d_k\|_2,\|G_ky_k\|_2\} &< \infty \label{eq.sum_x} \\
  \sum_{k = k_\epsilon}^\infty \alpha_k \max\{\|d_k\|_2^2,\|G_ky_k\|_2^2\} &< \infty. \label{eq.sum_alpha}
  \end{align}
  \esubequations
  In conjunction with \eqref{eq.eps_no_update}, the bound in \eqref{eq.sum_x} implies that the iterate sequence $\{x_k\}$ is a Cauchy sequence, meaning $\{x_k\} \to \xbar$ for some $\xbar \in \R{n}$.  At the same time, with \eqref{eq.eps_no_update}, the bound in \eqref{eq.sum_alpha} implies that $\{\alpha_k\} \to 0$.  We claim that this implies that $p_k = p$ for all sufficiently large $k \in \Kcal_\alpha$, where $\Kcal_\alpha$ is defined as in Lemma~\ref{lem.alpha}.  Indeed, since $\{\alpha_k\} \to 0$, it follows by Lemma~\ref{lem.well_posed}(c) that for sufficiently large $k \in \N{}$ either (i) $p_k < p$ and $\alpha_k = 0$ or (ii) $p_k = p$ and $\alpha_k > 0$.  Combined with the fact that $|\Kcal_d| < \infty$, it follows along with Lemma~\ref{lem.alpha} that there exists an infinite number of iterations indexed by $k \geq k_\epsilon$ such that $\alpha_k d_k \neq 0$ and $p_k = p$, whereas all other iterations for sufficiently large $k \geq k_\epsilon$ yield $\alpha_k = 0$.  Going forward, for ease of notation in the remainder of the proof of this case, since $x_{k+1} \gets x_k$ and $(H_{k+1},W_{k+1}) \gets (H_k,W_k)$ whenever $\alpha_k = 0$, let us proceed without loss of generality under the assumption that $k_\epsilon = 0$ and $\epsilon = \epsilon_0$, and that $\alpha_k > 0$, $d_k \neq 0$, and $p_k = p$ for all $k \in \N{}$.  Notice that under these conditions the set $\Kcal_{H,W}$ defined in Lemma~\ref{lem.well_posed} equals $\N{}$.  Correspondingly, for a given $\chi \in (0,1)$, let $\Kcal_\chi$ be the indices of iterations for which \eqref{eq.HW} holds; in particular, for $k \in \Kcal_\chi$, one has from \eqref{eq.HW_W1}--\eqref{eq.HW_W2} that
  \bequation\label{eq.mu}
  \max\{\|d_k\|_2^2,\|G_ky_k\|_2^2\} \leq \mu \|G_ky_k\|_{W_k}^2,\ \ \text{where}\ \ \mu := \max\left\{\frac{\overline\mu}{\underline\mu},\frac{1}{\underline\mu}\right\}.
  \eequation
  Since \eqref{eq.eps_update} does not hold for any $k \geq k_\epsilon$, it follows that either \eqref{eq.inexact_step_1} or \eqref{eq.inexact_step_2} holds for all $k \geq k_\epsilon$.  Hence, by Lemmas~\ref{lem.approx} and \ref{lem.approx_2}, it follows that \eqref{eq.theta_good} holds for all $k \geq k_\epsilon$, meaning for all $k \geq k_\epsilon$ that
  \bequation\label{eq.GyW_sigma}
  \|G_ky_k\|_{W_k} \leq (1 + \sigma_k)\|P_{W_k}(\conv(\{\nabla f(x)\}_{x\in\Xcal_k}))\|_{W_k}.
  \eequation
  \begin{hitemize}
  \item \textbf{Subcase 2a:} If $\xbar$ is $\epsilon$-stationary, then $\|P_{W_k}(\Gcal(\xbar,\epsilon_k))\|_{W_k} = 0$ for any $W_k \succ 0$.  Therefore, with $\mu \in (0,\infty)$ defined in \eqref{eq.mu}, $\omega = \nu \epsilon / (\sqrt{\mu}(1 + \sigma))$, and $(\zeta,\Tcal)$ chosen as in Lemma~\ref{lem.T}, it follows that there exists $k_\zeta \in \N{}$ with $k_\zeta \geq k_\epsilon$ such that $x_k \in \Bmbb(\xbar,\zeta)$ for all $k \geq k_\zeta$ and, with \eqref{eq.GyW_sigma},
  \bequation\label{eq.Gy_eps}
  \baligned
  \max\{\|d_k\|_2,\|G_ky_k\|_2\}
  &\leq \sqrt{\mu} \|G_ky_k\|_{W_k} \\
  &\leq \sqrt{\mu}(1 + \sigma_k)\|P_{W_k}(\conv(\{\nabla f(x)\}_{x\in\Xcal_k}))\|_{W_k} \\
  &\leq \sqrt{\mu}(1 + \sigma_k) \omega \leq \nu \epsilon
  \ealigned
  \eequation
  whenever $k \geq k_\zeta$, $k \in \Kcal_\chi$, and $\Xcal_k \in \Tcal$.  Combining \eqref{eq.eps_no_update} and \eqref{eq.Gy_eps}, it follows that $\Xcal_k \not\in \Tcal$ for all $k \geq k_\zeta$ with $k \in \Kcal_\chi$.  However, this is a probability zero event since for all such $k$ the set $\Xcal_k$ will contain new points from $\Bmbb(x_k,\epsilon_k)$ that are generated independently whether or not $k \in \Kcal_\chi$, meaning that with probability one there exists sufficiently large such $k$ with $k \in \Kcal_\chi$ and $\Xcal_k \in \Tcal$, which would yield \eqref{eq.Gy_eps}.
  \item \textbf{Subcase 2b:} If $\xbar$ is not $\epsilon$-stationary, then it follows from Lemma~\ref{lem.well_posed}(c) that $\alpha_k$ satisfies \eqref{eq.armijo} for all $k \in \N{}$.  In particular, \eqref{eq.armijo} holds either with $\alpha_k \geq \gamma \overline\alpha$ or with $\alpha_k < \gamma \overline\alpha$ such that
  \bequation\label{eq.bad_alpha}
  f(x_k + \gamma^{-1}\alpha_k d_k) - f(x_k) \geq - \underline\eta \gamma^{-1} \alpha_k \max\{\|d_k\|_2^2,\|G_ky_k\|_2^2\}.
  \eequation
  In the latter case, Lebourg's mean value theorem \cite[Theorem~2.3.7]{Clar83} implies the existence of a point $\xtilde_k \in [x_k,x_k + \gamma^{-1}\alpha_k d_k]$ and $\gtilde_k \in \eth f(\xtilde_k)$ such that
  \bequation\label{eq.lebourg}
  f(x_k + \gamma^{-1}\alpha_kd_k) - f(x_k) = \gamma^{-1} \alpha_k \gtilde_k^Td_k.
  \eequation
  Combining \eqref{eq.bad_alpha}, \eqref{eq.lebourg}, and the fact that $d_k = -W_kG_ky_k$, one finds that
  \bequation\label{eq.vWGy_1}
  \gtilde_k^TW_kG_ky_k \leq \underline\eta \max\{\|d_k\|_2^2,\|G_ky_k\|_2^2\}.
  \eequation
  On the other hand, for any $\omega \in (0,\infty)$ and $(\zeta,\Tcal)$ as in Lemma~\ref{lem.T}, there exists $k_\omega \geq k_\epsilon$ such that $x_k \in \Bmbb(\xbar,\min\{\zeta,\epsilon/3\})$ for $k \geq k_\omega$ and, with \eqref{eq.GyW_sigma},
  \bequation\label{eq.nonstationary}
  \baligned
  \|G_ky_k\|_{W_k}
  &\leq (1 + \sigma_k) \|P_{W_k}(\conv(\{\nabla f(x)\}_{x\in\Xcal_k}))\|_{W_k} \\
  &\leq (1 + \sigma_k) \|P_{W_k}(\Gcal(\xbar,\epsilon_k))\|_{W_k} + (1 + \sigma_k)\omega
  \ealigned
  \eequation
  whenever $k \geq k_\omega$, $k \in \Kcal_\chi$, and $\Xcal_k \in \Tcal$.  Hence, for such $k$, it follows by Lemma~\ref{lem.angle_lemma} with $\Scal = \Gcal(\xbar,\epsilon_k)$, $\beta = \underline\eta \mu \in (0,1)$ (where this inclusion is guaranteed by Algorithm~\ref{alg.sufficient_decrease}), and $W = W_k$ that for sufficiently small $\sigma_k \in (0,\sigma)$ and $\omega \in (0,\infty)$ one finds that \eqref{eq.mu} and \eqref{eq.nonstationary} imply
  \bequation\label{eq.vWGy_2}
  \baligned
    v^TW_kG_ky_k &> \underline\eta \mu \|G_ky_k\|_{W_k}^2 \\
    &\geq \underline\eta \max\{\|d_k\|_2^2,\|G_ky_k\|_2^2\}\ \ \text{for all}\ \ v \in \Gcal(\xbar,\epsilon_k).
  \ealigned
  \eequation
  There exists $k_\sigma \geq k_\omega$ such that $\sigma_k$ is sufficiently small for all $k \geq k_\sigma$ with $k \in \Kcal_\chi$ since the fact that $\{\alpha_k\} \to 0$ and the construction of the algorithm implies that $\{\sigma_k\} \to 0$.  Together, \eqref{eq.vWGy_1} and \eqref{eq.vWGy_2} imply that $\gtilde_k \notin \Gcal(\xbar,\epsilon_k)$ whenever $k \geq k_\sigma$, $k \in \Kcal_\chi$, and $\Xcal_k \in \Tcal$.  However, by the facts that $\mathds{1}^Ty_k = 1$ and $y_k \geq 0$, Assumption~\ref{ass.f}, and \cite[Proposition~2.1.2]{Clar83}, it follows for all $k \geq k_\sigma$ with $k \in \Kcal_\chi$ that
  \bequationNN
  \|d_k\|_2 = \|W_kG_ky_k\|_2 \leq \sqrt{\overline\mu} \|G_ky_k\|_2 \leq \sqrt{\overline\mu} L_{\Bmbb(\xbar,\epsilon)},
  \eequationNN
  where $L_{\Bmbb(\xbar,\epsilon)} \in (0,\infty)$ is a Lipschitz constant for $f$ over $\Bmbb(\xbar,\epsilon)$.  This shows that $\{\|d_k\|_2\}_{k \in \Kcal_\chi}$ is bounded.  This fact, along with $\{\alpha_k\} \to 0$, implies that $\alpha_k \leq \gamma \epsilon / (3 \|d_k\|_2)$ for all sufficiently large $k \in \Kcal_\chi$, i.e., $\gamma^{-1} \alpha_k \|d_k\|_2 \leq \epsilon/3$ for all sufficiently large $k \in \Kcal_\chi$.  Along with the fact that $x_k \in \Bmbb(\xbar,\min\{\zeta,\epsilon/3\})$ implies $\|x_k - \xbar\|_2 \leq \epsilon/3$, it follows that $\xtilde_k \in \Bmbb(\xbar,2\min\{\zeta,\epsilon/3\}/3)$ and hence $\gtilde_k \in \Gcal(\xbar,\epsilon_k)$ for all sufficiently large $k \in \N{}$.  Overall, since $\gtilde_k \not\in \Gcal(\xbar,\epsilon_k)$ whenever $k \geq k_\sigma$, $k \in \Kcal_\chi$, and $\Xcal_k \in \Tcal$, yet $\gtilde_k \in \Gcal(\xbar,\epsilon_k)$ for all sufficiently large $k$, it follows that $\Xcal_k \notin \Tcal$ for all sufficiently large $k \in \Kcal_\chi$.  However, this is a probability zero event since $|\Kcal_\chi| = \infty$ and the sample points are generated independently of whether $k \in \Kcal_\chi$.
  \end{hitemize}
  \end{hitemize}
  We have shown that $\{\epsilon_k\} \searrow 0$ with probability one.  If $\{\epsilon_k\} \searrow 0$, then by \eqref{eq.eps_update} there exists an infinite index set $\Kcal_\epsilon := \{k \in \N{} : \epsilon_{k+1} \gets \psi \epsilon_k\}$ where
  \bequationNN
  \max\{\|d_k\|_2,\|G_ky_k\|_2\} \leq \epsilon_k\ \ \text{for all}\ \ k \in \Kcal_\epsilon.
  \eequationNN
  The same argument as in \cite[Theorem~4.2,~Case~2]{CurtQue13}, which borrows from \cite[Theorem~3.3,~part~(iii)]{Kiwi07}, shows all cluster points of $\{x_k\}$ are stationary for $f$. \Halmos
\endproof

Our second convergence result, presented as the following corollary, considers the case when one chooses $\psi = 1$ so that the sampling radius remains that $\epsilon_0 \in (0,\infty)$ for all $k \in \Kcal$.  Similar results have appeared in the literature to prove a similar property of other GS methods; see, e.g., \cite[Theorem~3.5]{Kiwi07}.

\bcorollary\label{cor.gs_inexact}
  Suppose $\psi = 1$.  Algorithm~\ref{alg.gs_inexact} either terminates finitely with a stationary point for~$f$ or, with probability one, it performs an infinite number of outer iterations.  In the latter case, with probability one, it either reaches iteration $k \in \N{}$ such that $0 \in \Gcal(x_k,\epsilon_k)$ or every cluster point of the iterate sequence $\{x_k\}$ is $\epsilon_0$-stationary for~$f$.
\ecorollary
\proof{Proof.}
  As in the proof of Theorem~\ref{th.gs_inexact}, if Algorithm~\ref{alg.gs_inexact} terminates finitely with a stationary point for $f$, then there is nothing left to prove.  Otherwise, by Lemma~\ref{lem.well_posed}, it follows with probability one that an infinite number of outer iterations are performed, meaning $\Kcal = \N{}$.  If the algorithm reaches iteration $k \in \N{}$ in which $0 \in \Gcal(x_k,\epsilon_k)$, then there is nothing left to prove.  Otherwise, following the arguments in the proof of Theorem~\ref{th.gs_inexact}, it follows that $\inf \{\|G_ky_k\|_2 : k \in \N{}\} > 0$ is a probability zero event.  In the probability one event that $\inf \{\|G_ky_k\|_2 : k \in \N{}\} = 0$, the conclusion follows from the fact that $\eth_{\epsilon_0} f$ is closed. \Halmos
\endproof

\section{GS Algorithm with Gradient Aggregation}\label{sec.aggregation}

Our second algorithm adds a conceptually straightforward, but practically significant enhancement to Algorithm~\ref{alg.gs_inexact}.  In particular, we add a procedure for exploiting gradient aggregation that can significantly reduce the size of the subproblems to be solved approximately in each ``outer'' iteration of the algorithm.  We remark that this enhancement to the GS methodology is only possible when one is able to employ inexact subproblem solutions.  This is the case since the exact solution of a subproblem involving a ``gradient aggregation vector'' does not offer the exact solution of a subproblem involving individual gradients and no aggregation.

In this section, we present a statement of the proposed algorithm, then show that it offers the same convergence guarantees as does Algorithm~\ref{alg.gs_inexact}.

\subsection{Algorithm Description}

Our algorithm with inexact subproblem solutions and gradient aggregation is stated as Algorithm~\ref{alg.gs_aggregation}.  The algorithm borrows much from Algorithm~\ref{alg.gs_inexact}; we have written it in such a manner that only its unique steps are stated.  The main idea of the enhancement is the following.  For any $k + 1 \in \Kcal$ such that $\alpha_k > 0$, the matrix of gradients $G_{k+1}$ contains all points in the set $\Xcal_{k+1}$, as in Algorithm~\ref{alg.gs_inexact}.  However, for any $k+1 \in \Kcal$ such that $x_{k+1} = x_k$ since $\alpha_k = 0$, rather than solve a subproblem defined by gradients at all points in $\Xcal_{k+1}$, the algorithm considers a subproblem in which the gradients defining the matrix $G_k$ (which compose a submatrix of $G_{k+1}$) have been \emph{aggregated} into a single ``gradient aggregation vector'' $G_ky_k$.  The following lemma shows that a feasible point for the subproblem \emph{that the algorithm considers} in iteration $k+1$ corresponds to a feasible point for the subproblem \emph{that would be defined by all gradients in $G_{k+1}^{\rm full}$}.

\begin{algorithm}[ht]
  \caption{GS with Inexact Subproblem Solutions and Gradient Aggregation}
  \label{alg.gs_aggregation}
  \begin{algorithmic}[1]
  \Require [\dots same parameters and initial values as in Algorithm~\ref{alg.gs_inexact}, except $G_0$ \dots]
  \State Set $G_0^{\rm full}$ by \eqref{eq.G}, $G_0^{\rm agg}$ by \eqref{eq.G}, and $\alpha_{-1} \gets 0$.
  \For{\textbf{all} $k \in \N{}$}
  \If{$\alpha_{k-1} > 0$ or $p_k \geq p$}
  \State set $G_k \gets G_k^{\rm full}$;
  \Else
  \State set $G_k \gets G_k^{\rm agg}$.
  \EndIf
  \State [\dots same as Line~\ref{step.gs_inexact.terminate} through Line~\ref{step.HW} of Algorithm~\ref{alg.gs_inexact} \dots]
  \State Set $(\Xcal_{k+1},p_{k+1})$ by Algorithm~\ref{alg.sample} and $G_{k+1}^{\rm full}$ by \eqref{eq.G}.
  \If{$\alpha_k > 0$}
  \State set $G_{k+1}^{\rm agg} \gets G_{k+1}^{\rm full}$;
  \Else
  \State set $G_{k+1}^{\rm agg} \gets \bbmatrix \nabla f(x_{k+1}) & G_ky_k & [\nabla f(x)]_{x \in \Xcal_{k+1} \setminus (x_{k+1} \cup \Xcal_k)} \ebmatrix$.
  \EndIf
  \EndFor
  \end{algorithmic}
\end{algorithm}

\blemma\label{lem.feas}
  Consider $k \in \Kcal$ such that $k \geq 1$ and $\alpha_{k-1} = 0$, meaning $G_k = G_k^{\rm agg}$.  For any $j \in \N{}$ such that~$y_{k,j}$ is computed, this vector, which is feasible for the dual problem in \eqref{prob.qp}, corresponds uniquely to a feasible point for the dual problem in \eqref{prob.qp} if $G_k^{\rm full}$ were used in place of~$G_k = G_k^{\rm agg}$.
\elemma
\proof{Proof.}
  Consider any $j \in \N{}$ such that $y_{k,j}$ is computed.  Let $[y_{k,j}]_1$ and $[y_{k,j}]_2$ denote the first and second elements of $y_{k,j}$, respectively, with the subvector of all remaining elements of $y_{k,j}$ being denoted as $[y_{k,j}]_{>2}$.  One finds that
  \bequationNN
  \baligned
    &\ G_k^{\rm agg}y_{k,j} \\
   =&\ \nabla f(x_k)[y_{k,j}]_1 + (G_{k-1}y_{k-1})[y_{k,j}]_2 + [\nabla f(x)]_{x \in \Xcal_k \setminus (x_k \cup \Xcal_{k-1})} [y_{k,j}]_{>2} \\
   =&\ \bbmatrix \nabla f(x_k) & G_{k-1} & [\nabla f(x)]_{x \in \Xcal_k \setminus (x_k \cup \Xcal_{k-1})} \ebmatrix \bbmatrix [y_{k,j}]_1 \\ y_{k-1} [y_{k,j}]_2 \\ [y_{k,j}]_{>2} \ebmatrix = G_k^{\rm full} \bbmatrix [y_{k,j}]_1 \\ y_{k-1} [y_{k,j}]_2 \\ [y_{k,j}]_{>2} \ebmatrix,
  \ealigned
  \eequationNN
  where---since $\mathds{1}^Ty_{k-1} = 1$, $\mathds{1}^Ty_{k,j} = 1$, $y_{k-1} \geq 0$, and $y_{k,j} \geq 0$---it follows that
  \bequationNN
  \mathds{1}^T \bbmatrix [y_{k,j}]_1 \\ y_{k-1} [y_{k,j}]_2 \\ [y_{k,j}]_{>2} \ebmatrix = 1\ \ \text{and}\ \ \bbmatrix [y_{k,j}]_1 \\ y_{k-1} [y_{k,j}]_2 \\ [y_{k,j}]_{>2} \ebmatrix \geq 0,
  \eequationNN
  which proves the desired result. \Halmos
\endproof

\btheorem\label{th.gs_aggregation}
  Suppose $\psi \in (0,1)$.  Algorithm~\ref{alg.gs_aggregation} either terminates finitely with a stationary point for $f$ or, with probability one, it performs an infinite number of outer iterations.  In the latter case, with probability one, the sampling radius sequence satisfies $\{\epsilon_k\} \searrow 0$ and every cluster point of the iterate sequence $\{x_k\}$ is stationary for~$f$.
\etheorem
\proof{Proof.}
  For all $k \in \N{}$, the result of Lemma~\ref{lem.exact_solution} holds regardless of whether $G_k = G_k^{\rm agg}$ or $G_k = G_k^{\rm full}$ due to the fact that $G_k$ has $\nabla f(x_k)$ as its first column in either case.  The results of Lemmas~\ref{lem.approx} and~\ref{lem.approx_2} also continue to hold regardless of whether $G_k = G_k^{\rm agg}$ or $G_k = G_k^{\rm full}$, implying that the inner loop terminates finitely for all $k \in \Kcal$.  Now consider the pair $(d_k,y_k) = (-W_kG_ky_k,y_k)$ upon termination of the inner loop in iteration $k \in \Kcal$.  If $G_k = G_k^{\rm full}$, then the properties of $(d_k,y_k)$ are the same as that in Algorithm~\ref{alg.gs_inexact}.  Otherwise, when $G_k = G_k^{\rm agg}$, one may consider
  \bequation\label{eq.yk_agg}
  \bbmatrix [y_k]_1 \\ y_{k-1} [y_k]_2 \\ [y_k]_{>2} \ebmatrix
  \eequation
  as the dual vector, as shown by Lemma~\ref{lem.feas}.  The arguments of Lemmas~\ref{lem.well_posed}--\ref{lem.T} and Theorem~\ref{th.gs_inexact} now follow in the same manner as in Section~\ref{sec.inexact} using $G_k^{\rm full}$ in place of $G_k$ and $y_k$ or \eqref{eq.yk_agg} in place of the dual vector for all $k \in \Kcal$.  Crucial in these arguments is that, if the sample set size indicator $p_k$ ever exceeds $p$, then $G_k = G_k^{\rm full}$ and the algorithm behaves as Algorithm~\ref{alg.gs_inexact} for such $k \in \Kcal$. \Halmos
\endproof

\bcorollary\label{cor.gs_aggregation}
  Suppose $\psi = 1$.  Algorithm~\ref{alg.gs_aggregation} either terminates finitely with a stationary point for~$f$ or, with probability one, it performs an infinite number of outer iterations.  In the latter case, with probability one, it either reaches iteration $k \in \N{}$ such that $0 \in \Gcal(x_k,\epsilon_k)$ or every cluster point of the iterate sequence $\{x_k\}$ is $\epsilon_0$-stationary for~$f$.
\ecorollary
\proof{Proof.}
  The proof follows from that of Theorem~\ref{th.gs_aggregation} in the same manner as the proof of Corollary~\ref{cor.gs_inexact} follows from that of Theorem~\ref{th.gs_inexact}. \Halmos
\endproof

\section{Numerical Experiments}\label{sec.numerical}

In this section, we present the results of numerical experiments with implementations of our proposed algorithms.  The main purpose of these experiments is to show that the introduction of inexactness and gradient aggregation can reduce the computational expense of an adaptive GS algorithm consistently and often substantially.  As a sanity check, we also provide a comparison between our implementation of Algorithm~\ref{alg.gs_aggregation} and a state-of-the-art code.  All experiments were run on a Macbook Air with a 2.2 GHz Dual-Core Intel Core i7 processor running macOS 11.4.

We implemented our algorithms in the C++ software package \texttt{NonOpt} \citep{CurtNonOpt}.  For the parameters used in the algorithms and subroutines, we employed the values stated in Table~\ref{tab.parameters}.  These values are used consistently across all of our experiments.  As is typical in implementations of GS methods, our implementations assume that $x_k + \alpha_k d_k \in \Dcal$ for all $k \in \N{}$, meaning that the loop in Algorithm~\ref{alg.perturb} always terminates in the first iteration; hence, the parameter $\overline\ell$ is not used.  The initial point $x_0 \in \R{n}$ in each run of the algorithm was chosen in a problem-dependent manner; see the references given below in our discussion of the test problems used.

\btable[ht]
  \centering
  \caption{User-specified parameters for our implemented algorithms and subroutines.}
  \label{tab.parameters}
  \resizebox{\textwidth}{!}{%
  \btabular{|c|c|c|l|}
  \hline
  Parameter(s) & Range &Values & Description \\
  \hline
  $\nu$  & $(0,\infty)$ & $1$ & Stationarity measure tolerance \\
  $\underline\alpha \leq \overline\alpha$ & $(0,\infty)$& $10^{-20}\leq 100$ & Stepsize thresholds \\
  $\alpha_{\text{init}}$ & $(0,\infty)$ & $1$ & Initial stepsize \\
  $\rho$ & $(0,1)$ & $0.01$ & Inexactness threshold bound \\
  $\kappa$ & $(0,1)$ & $0.0001$ & Inexactness threshold \\
  $\psi$ & $(0,1)$ & $0.1$ & Sampling radius reduction factor \\
  $\iota$ & $(0,1)$ &$0.5$ & Inexactness parameter reduction factor \\
  $\underline\eta < \overline\eta$ & $(0,1)$ &$10^{-10} < 0.9$ & Armijo--Wolfe line search parameters \\
  $p$ & $[n+1,\infty)\cap\Nmbb$ & $10n$ & Sample set size threshold \\
  $\sigma$ & $(0,\infty)$ & $10$ & Inexactness threshold reset value \\
  $\gamma$ & $(0,1)$ & 0.5 & Stepsize modification factor \\
  $\underline\phi < 1 < \overline\phi$ & $(0,\infty)$ &$10^{-20} < 1 < 10^8$ & BFGS updating thresholds \\
  $\xi$ & $(0,\infty)$ & $10^{-20}$ & Curvature threshold \\
  $\overline p$ & $\N{}$ & 100 & Size of addition to sample set \\
  $H_0$ & $\succ 0$ & $I$ & Initial Hessian approximation \\
  $\epsilon_0$ & $(0,\infty)$ & $\max\{0.01,0.1\|\nabla f (x_0)\|_\infty\}$ & Initial stationarity radius \\
  \hline
  \etabular
  }
\etable

\texttt{NonOpt} contains a dual active-set QP solver that we used for solving the QP subproblems arising in the implementations of our algorithms.  To reduce CPU time, during the solve of a given QP, the termination conditions \eqref{eq.eps_update}--\eqref{eq.inexact_step_2} are not checked in every iteration of the QP solver.  Instead, these conditions are checked only after $(p_k+1)/4$ QP iterations have been performed, and after this threshold is reached, the conditions are checked only once every four QP solver iterations.

In our implementations, the outer iteration sequence terminates if
\bequation\label{eq.termination}
  \max\{\|G_ky_k\|_{\infty},\|W_kG_ky_k\|_{\infty},\epsilon_k\} \leq 10^{-4}
\eequation
or once an objective function value tolerance or CPU time limit is reached.  These latter criteria are discussed in further detail in the subsequent subsections.

We consider the performance of three implementations, to which we refer as follows:
\begin{hitemize}
  \item \texttt{GS-exact}: An implementation of an adaptive GS method in which the QP subproblems are solved ``exactly'' in each iteration; in particular, every aspect of this implementation is the same as that of \texttt{GS-inexact} (below), except that, when tasked to solve each QP subproblem, the QP solver is run until the $\ell_\infty$-norm of the KKT error for the QP is reduced below $10^{-10}$.
  \item \texttt{GS-inexact}: An implementation of Algorithm~\ref{alg.gs_inexact}.
  \item \texttt{GS-inexact-agg}: An implementation of Algorithm~\ref{alg.gs_aggregation}.
\end{hitemize}

\subsection{Randomly Generated Test Problems}

Our algorithms are designed to minimize objectives that may be nonconvex and/or nonsmooth.  However, in order to conduct a controlled comparison between the aforementioned implemented algorithms, our main experiment involves randomly generated convex test problems of the form
\bequationNN
  \min_{x \in \R{n}}\ g^Tx + \thalf x^THx + \max\{Ax + b\},
\eequationNN
where $g \in \R{n}$, $H \in \R{n \times n}$ is symmetric and positive definite, $A \in \R{m \times n}$, $b \in \R{m}$, and the $\max$ is taken element-wise.  (By employing convex, as opposed to nonconvex test problems, we can be sure that the results of our experiments are not skewed by two algorithms converging to different local minimizers and other related circumstances.)  The problems were constructed such that the unique global minimizer is always $x_* = 0$, the global minimum is always $f(x_*) = 0$, and the number of elements of the vector $Ax_* + b = b$ yielding the $\max$, call it $m_\Acal$, is always predetermined.

When solving potentially nonconvex and/or nonsmooth optimization problems, termination conditions can be sensitive in practice; e.g., one can find that the termination condition~\eqref{eq.termination} may be satisfied relatively early for some problems, whereas for other problems the magnification of small numerical errors can cause \eqref{eq.termination} to take longer to be satisfied.  Hence, we added a condition for these experiments that terminates an algorithm whenever the objective value is less than a prescribed threshold of $10^{-3}$.  This is reasonable in these experiments since $f(x_*) = 0$ for all problems.

For the purposes of these experiments, fifteen problems were generated; with $n = 1000$ and $m = 500$, five problems were generated for each of the values $m_\Acal \in \{125,250,375\}$.  In this manner, we provide results for a range of dimensions of the ``$\Ucal$-space'' and ``$\Vcal$-space'' at the minimizer; see, e.g., \citep{Liu2020}.  For each problem, each of the three implemented algorithms were run from the same randomly generated starting point; in particular, each element for the initial point was drawn from a standard normal distribution.  Since GS methods are randomized, we ran each algorithm 10 times for each problem and provide averages of performance measures over these 10 runs.

Results for \texttt{GS-exact}, \texttt{GS-inexact}, and \texttt{GS-inexact-agg} are provided in Tables~\ref{exact.results}, \ref{inexact.results}, and \ref{inexactagg.results}, respectively.  Averaged over the 10 runs for each algorithm and problem, we provide the required number of iterations (\texttt{iters}), required total number of QP solver iterations (\texttt{QP-iters}), required number of objective function evaluations (\texttt{funcs}), required number of objective gradient evaluations (\texttt{grads}), and final objective value ($f$).  Since the total computational effort is roughly proportional to the total number of QP iterations, for \texttt{GS-inexact} and \texttt{GS-inexact-agg}, we provide the relative change in the required total number of QP iterations as compared to \texttt{GS-exact}.  (This is a rough proxy for computational effort since the cost for each QP solver iteration can differ depending on the number of nonzero variables in the dual solution estimate.  That said, we found it to be the best measure for comparison, as opposed to CPU time which can vary despite the algorithm being run with the same initial conditions, random number generator seeds, and so on.)  In these statistics, a negative percentage indicates that \texttt{GS-inexact} (or \texttt{GS-inexact-agg}) required fewer total QP solver iterations than \texttt{GS-exact}; e.g., a statistic of $-z\%$ indicates that the algorithm lowered the required total number of QP solver iterations by $z\%$.

\begin{table}[ht]
  \centering
\caption{Results for \texttt{GS-exact} averaged over 10 runs.}
\label{exact.results}
  \footnotesize
\texttt{
\begin{tabular}{rrr|rrrrr}
$n$ & $m$ & $m_\Acal$ & iters & QP-iters & funcs & grads & $f$ \\
\hline
1000 & 500 &  125 &       964 &     4583 &     6486 &    21495 &  +9.948857e-04  \\
1000 & 500 &  125 &       961 &     5199 &     6521 &    24482 &  +9.943792e-04  \\
1000 & 500 &  125 &       898 &     3042 &     6023 &    19725 &  +9.956177e-04  \\
1000 & 500 &  125 &      1029 &     3026 &     6782 &    20148 &  +9.940536e-04  \\
1000 & 500 &  125 &       916 &     4040 &     6258 &    25183 &  +9.953463e-04  \\
1000 & 500 &  250 &      1076 &    16090 &     7433 &    41892 &  +9.922383e-04  \\
1000 & 500 &  250 &       869 &    10493 &     6303 &    36224 &  +9.893498e-04  \\
1000 & 500 &  250 &      1210 &    37689 &     7563 &    43487 &  +1.194211e-03  \\
1000 & 500 &  250 &      1088 &    13620 &     7189 &    38885 &  +1.409556e-03  \\
1000 & 500 &  250 &      1080 &    13870 &     7204 &    41252 &  +1.092632e-03  \\
1000 & 500 &  375 &      2063 &    47562 &    10405 &    55532 &  +2.080894e-03  \\
1000 & 500 &  375 &      1861 &    56658 &    10348 &    74029 &  +1.329274e-03  \\
1000 & 500 &  375 &      2193 &    87014 &    11563 &    74775 &  +1.756560e-03  \\
1000 & 500 &  375 &      2061 &    71832 &    11111 &    79926 &  +2.028260e-03  \\
1000 & 500 &  375 &      1882 &    65654 &    10059 &    62165 &  +2.326188e-03  \\
\end{tabular}
}
\end{table}

\begin{table}[ht]
  \centering
\caption{Results for \texttt{GS-inexact} averaged over 10 runs.  The final column indicates the relative change in \texttt{QP-iters} compared to \texttt{GS-exact}.}
\label{inexact.results}
  \footnotesize
\texttt{
\begin{tabular}{rrr|rrrrr|r}
$n$ & $m$ & $m_\Acal$ & iters & QP-iters & funcs & grads & $f$ & change in QP-iters \\
\hline
1000 & 500 &  125 &       957 &     3024 &     6518 &    26306 &  +9.968597e-04  & -34.003142\% \\
1000 & 500 &  125 &       948 &     3131 &     6448 &    26093 &  +9.943444e-04  & -39.767679\% \\
1000 & 500 &  125 &       897 &     2818 &     6116 &    23888 &  +9.951412e-04  & -7.350186\% \\
1000 & 500 &  125 &      1025 &     2848 &     6837 &    22942 &  +9.967565e-04  & -5.865829\% \\
1000 & 500 &  125 &       905 &     4096 &     6351 &    32989 &  +9.944478e-04  & +1.391089\% \\
1000 & 500 &  250 &      1105 &     9444 &     7456 &    44674 &  +1.078713e-03  & -41.305334\% \\
1000 & 500 &  250 &      1036 &     8875 &     7015 &    44131 &  +1.076791e-03  & -15.417969\% \\
1000 & 500 &  250 &       915 &     3953 &     6656 &    39961 &  +9.926499e-04  & -89.511723\% \\
1000 & 500 &  250 &       924 &     3428 &     6717 &    39259 &  +9.904635e-04  & -74.829662\% \\
1000 & 500 &  250 &      1082 &     8180 &     7559 &    48848 &  +1.031102e-03  & -41.021038\% \\
1000 & 500 &  375 &      1362 &    17336 &     8145 &    55212 &  +1.509360e-03  & -63.549472\% \\
1000 & 500 &  375 &      1547 &    31698 &     9817 &    83939 &  +1.379415e-03  & -44.053803\% \\
1000 & 500 &  375 &      1382 &    16480 &     8545 &    63765 &  +1.559232e-03  & -81.059730\% \\
1000 & 500 &  375 &      1718 &    35042 &     9584 &    66791 &  +1.977902e-03  & -51.216579\% \\
1000 & 500 &  375 &      2031 &    56392 &    11059 &    92092 &  +1.801977e-03  & -14.107727\% \\
\end{tabular}
}
\end{table}

\begin{table}[ht]
  \centering
\caption{Results for \texttt{GS-inexact-agg} averaged over 10 runs.  The final column indicates the relative change in \texttt{QP-iters} compared to \texttt{GS-exact}.}
\label{inexactagg.results}
  \footnotesize
\texttt{
\begin{tabular}{rrr|rrrrr|r}
$n$ & $m$ & $m_\Acal$ & iters & QP-iters & funcs & grads & $f$ & change in QP-iters \\
\hline
1000 & 500 &  125 &       974 &     2130 &     6428 &    17515 &  +9.950228e-04  & -53.521403\% \\
1000 & 500 &  125 &       965 &     2147 &     6413 &    18405 &  +9.965422e-04  & -58.709156\% \\
1000 & 500 &  125 &       911 &     1682 &     6020 &    14611 &  +9.948074e-04  & -44.679662\% \\
1000 & 500 &  125 &      1012 &     2715 &     6733 &    22958 &  +9.941496e-04  & -10.277594\% \\
1000 & 500 &  125 &       929 &     2022 &     6180 &    17479 &  +9.948088e-04  & -49.940594\% \\
1000 & 500 &  250 &       935 &     2775 &     6675 &    32963 &  +9.948877e-04  & -82.748634\% \\
1000 & 500 &  250 &      1021 &    10064 &     6742 &    33816 &  +1.002834e-03  & -4.092971\% \\
1000 & 500 &  250 &       897 &     2316 &     6350 &    30019 &  +9.920159e-04  & -93.854024\% \\
1000 & 500 &  250 &       932 &     2541 &     6608 &    31248 &  +9.937333e-04  & -81.339941\% \\
1000 & 500 &  250 &       924 &     2302 &     6459 &    27277 &  +9.965494e-04  & -83.398221\% \\
1000 & 500 &  375 &       860 &     2679 &     6414 &    45224 &  +9.910356e-04  & -94.367142\% \\
1000 & 500 &  375 &       855 &     2702 &     6466 &    47196 &  +9.928050e-04  & -95.229876\% \\
1000 & 500 &  375 &       831 &     2676 &     6285 &    47392 &  +9.938719e-04  & -96.924309\% \\
1000 & 500 &  375 &       960 &    14692 &     6797 &    48239 &  +1.141195e-03  & -79.545916\% \\
1000 & 500 &  375 &      1015 &    19647 &     7155 &    54267 &  +1.005219e-03  & -70.073932\% \\
\end{tabular}
}
\end{table}

Let us start with a few general observations about these experiments.  First, between the termination condition \eqref{eq.termination} and the condition that the algorithm terminates if the objective value fell below $10^{-3}$, one finds that the solutions obtained by all algorithms on all problems were comparable in quality with final objective values on the order of $10^{-3}$.  Second, one finds that the computational effort required by the algorithms was directly proportional to $m_{\Acal}$; this is expected since increasing $m_{\Acal}$ increases the dimension of the $\Vcal$-space of the objective function at the solution.

Most importantly, one finds in the results in Tables~\ref{exact.results}, \ref{inexact.results}, and \ref{inexactagg.results} that inexactness and gradient aggregation reduces the total number of QP solver iterations consistently and often substantially.  Interestingly, one also finds in many cases that \texttt{GS-inexact} and \texttt{GS-inexact-agg} also require fewer outer iterations.  This was not necessarily expected, and might not represent behavior that one should anticipate in general.  That said, one explanation for this behavior is that requiring exact subproblem solutions may tend to produce shorter search directions, whereas by allowing inexactness in the subproblem solutions the algorithm is able to take longer steps in each iteration.  In any case, due to the reduced number of QP solver iterations required per outer iteration, one may expect a reduction in total computational effort for \texttt{GS-inexact} and \texttt{GS-inexact-agg} even if these algorithms were to require the same number, or even more, outer iterations than \texttt{GS-exact}.

\subsection{Test Set Problems}

To demonstrate that our implementations can be competitive with a state-of-the-art solver, we performed an experiment to compare the performance of the state-of-the-art code \texttt{LMBM} \citep{KarmLMBM} and \texttt{GS-inexact-agg}.  The experiments with \texttt{GS-inexact-agg} in the previous subsection were performed with full BFGS approximations, but the experiments in this subsection were performed with a limited-memory BFGS strategy with a history of 50 so that the algorithm would be more similar to \texttt{LMBM}, which uses limited memory approximations (with a history of 7).

We chose a set of 20 test problems for which \texttt{LMBM} has been tuned, some of which are convex and some of which are nonconvex.  The first ten problems come from \citep{HaarMietMaek04} and the second ten come from \citep{LuksTumaSiskVlceRame02}.  (\texttt{LMBM} comes with implementations of the first ten problems; for the remaining test problems, we obtained Fortran implementations from~\cite{LuksTestProblems}.)  In these sources, each problem is provided with an initial point $x_0 \in \R{n}$, which were the initial points that we used in our experiments.  All of the problems are scalable in the sense that they are defined for any value of $n \in \N{}$.  We chose $n = 1000$ for all problems.

\texttt{LMBM} and \texttt{GS-inexact-agg} have many differences.  For example, \texttt{LMBM} employs a bundle method while \texttt{GS-inexact-agg} employs a GS method.  The termination criteria of the two codes are also very different; e.g., besides observing termination criteria related to detecting stationarity, \texttt{LMBM} may terminate due to various reasons related to the iterate and/or objective value not changing significantly between iterations.  Hence, in order to offer a fair and illustrative comparison, we ran \texttt{LMBM} for each problem 10 times and found the average CPU time required per problem.  (The solve for most problems terminated within a few seconds; the only exceptions were that for \texttt{Test29\_5}, which required approximately 30 seconds, and that for \texttt{Test29\_13}, which required approximately 15 seconds.)  We then ran \texttt{GS-inexact-agg} with a CPU time limit of the maximum of the average time required by \texttt{LMBM} and one second.  This caused \texttt{GS-inexact-agg} to terminate in a large majority of the runs due to the CPU time limit, even though it would have continued to iterate to obtain better solutions if it were allowed to do so.  Due to the random behavior of \texttt{GS-inexact-agg}, we ran the solver 10 times for each problem and present averages over these runs.

The results obtained by the codes are shown in Table~\ref{tab.lmbm}.  \texttt{LMBM} reports the number of iterations (\texttt{iters}), function evaluations (\texttt{funcs}), and final objective value ($f$).  For \texttt{GS-inexact-agg}, we additionally provide (averaged over all runs) the number of QP subproblem solver iterations (\texttt{QP-iters}) and gradient evaluations (\texttt{grads}).  Since the algorithms have various differences, it is not necessarily informative to compare the number of iterations or function evaluations required by the two methods.  On the other hand, one can compare final objective values, with respect to which one finds that the results are generally comparable.  \texttt{LMBM} yields lower values for some problems while \texttt{GS-inexact-agg} yields lower values for a few others.

\begin{table}[ht]
  \centering
  \caption{Results for \texttt{LMBM} and \texttt{GS-inexact-agg} averaged over 10 runs}
  \label{tab.lmbm}
  \footnotesize
  \texttt{
  \btabular{l|rrr|rrrrr}
    \hline
    name & iters & funcs & \multicolumn{1}{c|}{$f$} & iters & QP-iters & funcs & grads & \multicolumn{1}{c}{$f$} \\
\hline
                MaxQ &    21940 &    22808 &  +4.987830e-06 &     5234 &     5253 &    15000 &     5267 &  +4.021445e-04 \\
              MxHilb &      441 &      861 &  +6.166410e-03 &      140 &      140 &      687 &      200 &  +1.134878e-03 \\
         Chained\_LQ &      300 &     1824 &  -1.412780e+03 &       66 &      150 &      446 &      130 &  -1.412639e+03 \\
     Chained\_CB3\_1 &      291 &     1690 &  +1.998000e+03 &       72 &      143 &      501 &      145 &  +2.012761e+03 \\
     Chained\_CB3\_2 &       66 &      150 &  +1.998000e+03 &       61 &       80 &      300 &       93 &  +1.998000e+03 \\
         ActiveFaces &      523 &      569 &  +1.376680e-14 &       15 &       17 &      391 &      344 &  +3.961526e-05 \\
  Brown\_Function\_2 &      493 &     4217 &  +2.136910e-09 &       27 &      141 &      201 &      107 &  +9.041635e-01 \\
 Chained\_Mifflin\_2 &      546 &     3892 &  -7.064510e+02 &       40 &      136 &      316 &      168 &  -7.062611e+02 \\
Chained\_Crescent\_1 &      177 &      817 &  +3.681010e-08 &       37 &       44 &      187 &       52 &  +4.728769e-08 \\
Chained\_Crescent\_2 &      903 &     9626 &  +1.369240e-04 &       34 &      141 &      262 &      101 &  +1.353061e-01 \\
           Test29\_2 &       62 &       63 &  +9.815390e-01 &      307 &      354 &     1643 &      358 &  +7.212000e-01 \\
           Test29\_5 &     1230 &     4563 &  +6.434430e-06 &      185 &      369 &      948 &      374 &  +8.359626e-07 \\
           Test29\_6 &       44 &       48 &  +2.000000e+00 &       37 &      142 &      292 &      119 &  +2.007236e+00 \\
          Test29\_11 &      283 &     1336 &  +1.203580e+04 &       13 &      147 &      170 &      121 &  +1.208292e+04 \\
          Test29\_13 &     3747 &     7092 &  +5.665460e+02 &       92 &     1451 &     1219 &     1134 &  +6.140429e+02 \\
          Test29\_17 &      962 &     2247 &  +3.574260e-03 &        8 &      152 &      148 &      114 &  +1.096233e-03 \\
          Test29\_19 &      143 &     1012 &  +1.000000e+00 &       29 &      144 &      274 &      123 &  +1.013406e+00 \\
          Test29\_20 &      277 &     3087 &  +5.000010e-01 &       46 &      140 &      394 &      142 &  +5.008065e-01 \\
          Test29\_22 &       21 &      172 &  +1.966970e-06 &       10 &      153 &      163 &      118 &  +3.056232e-04 \\
          Test29\_24 &      315 &     1945 &  +4.232150e-02 &       32 &      138 &      352 &      119 &  +1.102612e-01 \\
\hline
  \etabular
  }
\end{table}

\section{Conclusion}\label{sec.conclusion}

We have proposed, analyzed, and tested two algorithms for minimizing locally Lipschitz objective functions.  The algorithms are based on the gradient sampling methodology.  The unique feature of the first algorithm is that it can allow \emph{inexactness} in the subproblem solutions while maintaining convergence guarantees, which is new to the literature on gradient sampling methods.  The unique feature of the second algorithm is that it can use inexact subproblem solutions and \emph{aggregated} gradients in place of individual gradients in the subproblem definitions.  Our numerical experiments show that employing inexactness and aggregation can each reduce computational effort.

\bibliographystyle{plain}
\bibliography{references}

\begin{thebibliography}{10}

\bibitem{ApkaNollProt08}
P.~Apkarian, D.~Noll, and O.~Prot.
\newblock A trust region spectral bundle method for nonconvex eigenvalue
  optimization.
\newblock {\em SIAM Journal on Optimization}, 19(1):281--306, 2008.

\bibitem{Bert09}
D.~P. Bertsekas.
\newblock {\em Convex Optimization Theory}.
\newblock Athena Scientific, Nashua, NH, USA, 2009.

\bibitem{BurkLewiOver02}
J.~V. Burke, A.~S. Lewis, and M.~L. Overton.
\newblock {Approximating Subdifferentials by Random Sampling of Gradients}.
\newblock {\em Mathematics of Operations Research}, 27(3):567--584, 2002.

\bibitem{BurkLewiOver05}
J.~V. Burke, A.~S. Lewis, and M.~L. Overton.
\newblock {A Robust Gradient Sampling Algorithm for Nonsmooth, Nonconvex
  Optimization}.
\newblock {\em SIAM Journal on Optimization}, 15(3):751--779, 2005.

\bibitem{BurkCurtLewiOverSimo19}
James~V. Burke, Frank~E. Curtis, Adrian~S. Lewis, Michael~L. Overton, and Lucas
  E. A.~Sim\ {o}es.
\newblock {Gradient Sampling Methods for Nonsmooth Optimization}.
\newblock In {\em {Numerical Nonsmooth Optimization}}, chapter~6, pages
  201--225. Springer, 2020.

\bibitem{ByrdNoce89}
R.~H. Byrd and J.~Nocedal.
\newblock A tool for the analysis of quasi-{Newton} methods with application to
  unconstrained minimization.
\newblock {\em SIAM Journal on Numerical Analysis}, 26(3):727--739, 1989.

\bibitem{Clar83}
F.~H. Clarke.
\newblock {\em {Optimization and Nonsmooth Analysis}}.
\newblock Canadian Mathematical Society Series of Monographs and Advanced
  Texts. John Wiley \& Sons, New York, NY, USA, 1983.

\bibitem{CurtNonOpt}
F.~E. Curtis.
\newblock {NonOpt}.
\newblock \url{https://coral.ise.lehigh.edu/frankecurtis/nonopt/}, 2021.

\bibitem{CurtOver12}
F.~E. Curtis and M.~L. Overton.
\newblock {A Sequential Quadratic Programming Algorithm for Nonconvex,
  Nonsmooth Constrained Optimization}.
\newblock {\em {SIAM Journal on Optimization}}, 22(2):474--500, 2012.

\bibitem{CurtQue13}
Frank~E. Curtis and Xiaocun Que.
\newblock {An Adaptive Gradient Sampling Algorithm for Nonsmooth Optimization}.
\newblock {\em {Optimization Methods and Software}}, 28(6):1302--1324, 2013.

\bibitem{CurtQue15}
Frank~E. Curtis and Xiaocun Que.
\newblock {A Quasi-Newton Algorithm for Nonconvex, Nonsmooth Optimization with
  Global Convergence Guarantees}.
\newblock {\em {Mathematical Programming Computation}}, 7(4):399--428, 2015.

\bibitem{CurtRobiZhou19}
Frank~E. Curtis, Daniel~P. Robinson, and Baoyu Zhou.
\newblock {A Self-Correcting Variable-Metric Algorithm Framework for Nonsmooth
  Optimization}.
\newblock {\em {IMA Journal of Numerical Analysis}}, 40(2):1154--1187, 2020.

\bibitem{Gold77}
A.~A. Goldstein.
\newblock Optimization of {Lipschitz} continuous functions.
\newblock {\em Mathematical Programming}, 13(1):14--22, 1977.

\bibitem{HaarMietMaek04}
N.~Haarala, K.~Miettinen, and M.~M. M{\"a}kel{\"a}.
\newblock New limited memory bundle method for large-scale nonsmooth
  optimization.
\newblock {\em Optimization Methods and Software}, 19(6):673--692, 2004.

\bibitem{HaarMietMaek07}
N.~Haarala, K.~Miettinen, and M.~M. M{\"a}kel{\"a}.
\newblock Globally convergent limited memory bundle method for large-scale
  nonsmooth optimization.
\newblock {\em Mathematical Programming}, 109(1):181--205, 2007.

\bibitem{HareSaga10}
W.~Hare and C.~Sagastiz\'abal.
\newblock A redistributed proximal bundle method for nonconvex optimization.
\newblock {\em SIAM Journal on Optimization}, 20(5):2442--2473, 2010.

\bibitem{HeloSantSimo17}
E.~S. Helou, S.~A. Santos, and L.~E.~A. Sim{\~o}es.
\newblock {On the Local Convergence Analysis of the Gradient Sampling Method
  for Finite Max-Functions}.
\newblock {\em Journal of Optimization Theory and Applications},
  175(1):137--157, 2017.

\bibitem{HiriLema93b}
J.-B. Hiriart-Urruty and C.~Lemar\'{e}chal.
\newblock {\em {Convex Analysis and Minimization Algorithms II}}.
\newblock A Series of Comprehensive Studies in Mathematics. Springer-Verlag,
  New York, NY, USA, 1993.

\bibitem{HossUsch17}
S.~Hosseini and A.~Uschmajew.
\newblock {A Riemannian Gradient Sampling Algorithm for Nonsmooth Optimization
  on Manifolds}.
\newblock {\em SIAM Journal on Optimization}, 27(1):173--189, 2017.

\bibitem{KarmLMBM}
N.~Karmitsa.
\newblock {LMBM}.
\newblock \url{http://napsu.karmitsa.fi/lmbm}, accessed 2021.

\bibitem{Kiwi85}
K.~C. Kiwiel.
\newblock A linearization algorithm for nonsmooth minimization.
\newblock {\em Mathematics of Operations Research}, 10(2):185--194, 1985.

\bibitem{Kiwi85b}
K.~C. Kiwiel.
\newblock {\em Methods of Descent for Nondifferentiable Optimization}.
\newblock Lecture Notes in Mathematics. Springer-Verlag, New York, NY, USA,
  1985.

\bibitem{Kiwi96}
K.~C. Kiwiel.
\newblock Restricted step and {Levenberg-Marquardt} techniques in proximal
  bundle methods for nonconvex nondifferentiable optimization.
\newblock {\em SIAM Journal on Optimization}, 6(1):227--249, 1996.

\bibitem{Kiwi07}
K.~C. Kiwiel.
\newblock {Convergence of the Gradient Sampling Algorithm for Nonsmooth
  Nonconvex Optimization}.
\newblock {\em SIAM Journal on Optimization}, 18(2):379--388, 2007.

\bibitem{LemaNemiNest95}
C.~Lemar{\'e}chal, A.~Nemirovskii, and {Yu.} Nesterov.
\newblock New variants of bundle methods.
\newblock {\em Mathematical Programming}, 69(1):111--147, 1995.

\bibitem{LewiOver13}
A.~S. Lewis and M.~L. Overton.
\newblock {Nonsmooth Optimization via Quasi-Newton Methods}.
\newblock {\em Mathematical Programming}, 141(1--2):135--163, 2013.

\bibitem{Liu2020}
Shuai Liu and Claudia Sagastiz{\'a}bal.
\newblock {\em Beyond First Order: {$\Vcal\Ucal$}-Decomposition Methods}, pages
  297--329.
\newblock Springer International Publishing, Cham, 2020.

\bibitem{LuksTestProblems}
L.~Luksan.
\newblock {Test Problems in Fortran}.
\newblock \url{http://www.cs.cas.cz/~luksan/test.html}, accessed 2021.

\bibitem{LuksVlce98}
L.~Luk{\v{s}}an and J.~Vl{\v{c}}ek.
\newblock A bundle-{Newton} method for nonsmooth unconstrained minimization.
\newblock {\em Mathematical Programming}, 83(1):373--391, 1998.

\bibitem{LuksTumaSiskVlceRame02}
L.~Luk\v{s}an, M.~T\.uma, M.~\v{S}i\v{s}ka, J.~Vl\v{c}ek, and
  N.~Rame\v{s}ov\'a.
\newblock {UFO 2002: Interactive System for Universal Functional Optimization}.
\newblock Technical Report 883, Institute of Computer Science, Academy of
  Sciences of the Czech Republic, 2002.

\bibitem{MortMost19}
Morteza Maleknia and Mostafa Shamsi.
\newblock A gradient sampling method based on ideal direction for solving
  nonsmooth optimization problems.
\newblock {\em Journal of Optimization Theory and Applications},
  187(3):181--204, 2020.

\bibitem{Miff77}
R.~Mifflin.
\newblock An algorithm for constrained optimization with semismooth functions.
\newblock {\em Mathematics of Operations Research}, 2(2):191--207, 1977.

\bibitem{Miff82}
R.~Mifflin.
\newblock A modification and an extension of {Lemarechal's} algorithm for
  nonsmooth minimization.
\newblock In D.~C. Sorensen and R.~J.-B. Wets, editors, {\em Nondifferential
  and Variational Techniques in Optimization}, pages 77--90. Springer Berlin
  Heidelberg, Berlin, Heidelberg, 1982.

\bibitem{MiffSaga05}
R.~Mifflin and C.~Sagastiz\'abal.
\newblock A {$\Vcal\Ucal$}-algorithm for convex minimization.
\newblock {\em Mathematical Programming}, 104(2):583--608, 2005.

\bibitem{NoceWrig06}
J.~Nocedal and S.J. Wright.
\newblock {\em Numerical Optimization}.
\newblock Springer Series in Operations Research. Springer, New York, NY, USA,
  2nd edition, 2006.

\bibitem{Rusz06}
A.~Ruszczynski.
\newblock {\em {Nonlinear Optimization}}.
\newblock Princeton University Press, Princeton, NJ, USA, 2006.

\bibitem{SchrZowe92}
H.~Schramm and J.~Zowe.
\newblock A version of the bundle idea for minimizing a nonsmooth function:
  {Conceptual} idea, convergence analysis, numerical results.
\newblock {\em SIAM Journal on Optimization}, 2(1):121--152, 1992.

\bibitem{TangLiuJianLi14}
C.-M. Tang, S.~Liu, J.-B. Jian, and J.-L. Li.
\newblock {A Feasible SQP-GS Algorithm for Nonconvex, Nonsmooth Constrained
  Optimization}.
\newblock {\em Numerical Algorithms}, 65(1):1--22, Jan 2014.

\end{thebibliography}

\appendix
\normalsize

\section{Online Companion: GS Algorithm Subroutines}\label{app.subroutines}

In this online companion, we present subroutines needed for the GS algorithms with inexact subproblem solutions and gradient aggregation proposed in the paper.  Besides the first, these subroutines have been motivated and presented in previous articles, as mentioned in the following subsections.  We include them here, along with pointers to lemmas that prove their properties, for ease of reference for the main body of the paper.

\subsection{Sufficient decrease parameter selection}\label{sec.sufficient_decrease}

For our proposed algorithm, the sufficient decrease parameter $\underline\eta \in (0,1)$ employed in the Armijo-Wolfe conditions, specifically in \eqref{eq.armijo}, needs to be set sufficiently small relative to an upper bound on the condition number of (at least a subset of) the generated Hessian approximations.  A similar relationship was required in \cite[eq.~(3.1)]{CurtQue15}.  Generally speaking, the Hessian approximations in a BFGS updating strategy can become arbitrarily ill-conditioned, but it is sufficient for our purposes to employ an upper bound that holds for a fraction of \emph{good} iterations that will provably occur; see \cite{ByrdNoce89}.

\balgorithm[ht]
  \caption{Sufficient Decrease Parameter Selection}
  \label{alg.sufficient_decrease}
  \balgorithmic[1]
    \Require $\underline\phi \in (0,1)$, $\overline\phi \in (1,\infty)$, and $H_0 \succ 0$ from outer algorithm; $\chi \in (0,1)$.
    \State Set
    \bequationNN
      c_0 \gets \frac{1}{1 - \chi}\(\trace(H_0) - \ln\det(H_0) + \overline\phi - 1 - \ln\underline\phi\) \in (0,\infty).
    \eequationNN
    \State Set $c_1 \gets e^{-c_0/2}$.
    \State Set $c_2$ (resp.,~$c_3$) as the smallest (resp.,~largest) value in $(0,\infty)$ such that
    \bequationNN
      1 - r + \ln r \geq -c_0\ \ \text{for all}\ \ r \in [c_2,c_3].
    \eequationNN
    \State Set $\underline\mu \gets c_1^2/c_3$ and $\overline\mu \gets 1/c_2^2$.
    \State Set
    \bequationNN
      \mu \gets \max\left\{\frac{\overline\mu}{\underline\mu},\frac{1}{\underline\mu}\right\} \in (1,\infty).
    \eequationNN
    \State Choose $\underline\eta \in (0,1)$ such that $\underline\eta \mu \in (0,1)$.
    \State \textbf{terminate} and \textbf{return} $(\underline\eta,\mu)$
  \ealgorithmic
\ealgorithm

\subsection{Line search}

Given a descent direction for $f$ at $x_k$, the line search is intended to find a stepsize satisfying the weak Armijo-Wolfe conditions; see \eqref{eq.wolfe} and, e.g., \cite{NoceWrig06}.  However, as motivated by \cite{CurtQue15}, the line search may terminate early if the sample set size indicator $p_k$ is less than the prescribed integer $p \in [n+1,\infty)$ or may ignore the curvature condition \eqref{eq.curv}---and switch to a ``backtracking Armijo'' line search---if both $p_k \geq p$ and a certain number of iterations of the line search have already been performed without termination.  This potential switch to a backtracking Armijo line search may be needed due to potential nonsmoothness of~$f$, since under Assumption~\ref{ass.f} one can only show that an Armijo-Wolfe line search can ``bracket'' a stepsize satisfying the Armijo-Wolfe conditions \eqref{eq.wolfe}; see \cite[Theorem~4.7]{LewiOver13}.  One could ensure finite termination of an Armijo-Wolfe line search, without having to switch to a backtracking Armijo line search as a backup, with a stronger assumption on $f$, such as it being weakly lower semismooth; see \citep{Miff77}.

Our line search subroutine is stated as Algorithm~\ref{alg.ls}.  Given $\underline\eta \in (0,1)$ and $\overline\eta \in (\underline\eta,1)$, the Armijo-Wolfe conditions that we use are stated as \eqref{eq.wolfe}.  We remark that \eqref{eq.armijo} does not use the directional derivative of $f$ at $x_k$ along~$d_k$, as is typical in the Armijo condition in the context of smooth optimization; rather, it uses squared norms of search direction quantities, which is common in GS algorithms when nonnormalized search directions are used; see, e.g., \cite[Eq.~(4.2)]{Kiwi07}.

\balgorithm[ht]
  \caption{Armijo-Wolfe Line Search}
  \label{alg.ls}
  \balgorithmic[1]
    \Require $(\underline\eta,\overline\eta,\underline\alpha)$ from outer algorithm; $\overline\alpha \in [\underline\alpha,\infty)$; $\alpha_{\text{init}} \in (0,\infty)$; $\gamma \in (0,1)$.
    \State Set $l \gets 0$, $u \gets \overline\alpha$, and $\alpha_k \gets \alpha_{\text{init}}$.
    \If{$d_k=0$}
      \State \textbf{terminate} and \textbf{return} $\alpha_k$.
    \EndIf
    \For{\textbf{all} $\ell \in \N{}$}
      \If{$p_k < p$ and $\alpha_k < \underline\alpha$}
        \State set $\alpha_k \gets 0$, \textbf{terminate}, and \textbf{return} $\alpha_k$. \Comment{truncate and take null stepsize}
      \EndIf
      \If{$\alpha_k < \underline\alpha$}
        \State set $l \gets 0$. \Comment{switch to backtracking Armijo line search}
      \EndIf
      \If{\eqref{eq.wolfe} holds or both $\alpha_k < \underline\alpha$ and \eqref{eq.armijo} hold}
        \State \textbf{terminate} and \textbf{return} $\alpha_k$. \Comment{success}
      \EndIf
      \If{\eqref{eq.armijo} does not hold}
        \State set $u \gets \alpha_k$;
      \Else
        \State set $l \gets \alpha_k$.
      \EndIf
      \State set $\alpha_k \gets (1 - \gamma)l + \gamma u$.
    \EndFor
  \ealgorithmic
\ealgorithm

A proof of Lemma~\ref{lem.well_posed}(c) can be found in that for \cite[Lemma~2.3]{CurtQue15}.

\subsection{Iterate perturbation}

Gradient sampling algorithms require that each iterate lies in the set of points over which the objective function $f$ is differentiable (for well-posedness of the algorithm) or even continuously differentiable (for the convergence guarantees); see \cite{BurkCurtLewiOverSimo19}.  For our proposed algorithms, we employ Algorithm~\ref{alg.perturb} to ensure that each iterate lies in the set~$\Dcal$ defined in Assumption~\ref{ass.f}.  If, after the line search, the resulting trial point is contained in $\Dcal$, i.e., $x_k + \alpha_k d_k \in \Dcal$, then $x_{k+1}$ is set to be this trial point; otherwise, the iterate perturbation strategy in Algorithm~\ref{alg.perturb} aims to compute $x_{k+1} \in \Dcal$ satisfying \eqref{eq.wolfe2}.

\balgorithm[ht]
  \caption{Iterate Perturbation}
  \label{alg.perturb}
  \balgorithmic[1]
    \Require $\overline\ell \in \N{}$.
    \State Set $x_{k+1} \gets x_k + \alpha_k d_k$.
    \If{$\alpha_k = 0$ or $d_k = 0$}
      \State \textbf{terminate} and \textbf{return} $x_{k+1}$.
    \EndIf
    \For{\textbf{all} $\ell \in \N{}$}
      \If{$x_{k+1} \in \Dcal$ and either \eqref{eq.wolfe2} holds or each of \eqref{eq.armijo2}, \eqref{eq.perturb}, and $\ell > \overline\ell$ hold}
        \State \textbf{terminate} and \textbf{return} $x_{k+1}$
      \EndIf
      \State Sample $x_{k+1}$ from a uniform distribution over
      \bequationNN
        \Bmbb\(x_k + \alpha_kd_k,\frac{\min\{\alpha_k,\epsilon_k\}\min\{\|d_k\|_2,\|G_ky_k\|_2\}}{\ell\max\{\|d_k\|_2,\|G_ky_k\|_2\}}\).
      \eequationNN
    \EndFor
  \ealgorithmic
\ealgorithm

Algorithm~\ref{alg.perturb} can fail if its \textbf{for} loop iterates infinitely.  However, under Assumption~\ref{ass.f}, this is a probability zero event.  In other words, it terminates finitely---meaning the subroutine runs successfully---with probability one.

A proof of Lemma~\ref{lem.well_posed}(d) can be found in that for \cite[Lemma~2.3]{CurtQue15}.

\subsection{Hessian and inverse Hessian approximation strategy}\label{sec.hessian}

The Hessian approximation strategy employed in \cite{CurtQue15} is conservative in the sense that it might replace a BFGS approximation with an L-BFGS approximation in order to ensure that, in certain cases, the eigenvalues of the Hessian approximation are bounded above and below away from zero.  For our purposes, we employ the less conservative strategy advocated in \cite{CurtRobiZhou19}, which exploits the \emph{self-correcting} properties of BFGS updating.  The subroutine we use is stated in Algorithm~\ref{alg.hessian}.

\balgorithm[ht]
  \caption{Hessian and Inverse Hessian Approximation Updates}
  \label{alg.hessian}
  \balgorithmic[1]
    \Require $\underline\phi \in (0,1)$ and $\overline\phi \in (1,\infty)$ from outer algorithm.
    \State Set $s_k \gets \alpha_kd_k$ and $\vhat_k \gets x_{k+1} - x_k$.
    \If{$s_k = 0$}
      \State set $(H_{k+1},W_{k+1}) \gets (H_k,W_k)$.
    \Else
      \State Compute $\vartheta_k$ as the smallest value in $[0,1]$ such that
      \bequationNN
        v_k \gets \vartheta_k s_k + (1 - \vartheta_k) y_k
      \eequationNN
      \State yields
      \bequation\label{eq.phi}
        \underline\phi \leq \frac{s_k^Tv_k}{\|s_k\|_2^2}\ \ \text{and}\ \ \frac{\|v_k\|_2^2}{s_k^Tv_k} \leq \overline\phi,
      \eequation
      \State then set
      \bequationNN
        \baligned
          H_{k+1} &\gets \(I - \frac{s_ks_k^TH_k}{s_k^TH_ks_k}\)^TH_k\(I - \frac{s_ks_k^TH_k}{s_k^TH_ks_k}\) + \frac{v_kv_k^T}{s_k^Tv_k} \\ \text{and}\ \ 
          W_{k+1} &\gets \(I - \frac{v_ks_k^T}{s_k^Tv_k}\)^TW_k\(I - \frac{v_ks_k^T}{s_k^Tv_k}\) + \frac{s_ks_k^T}{s_k^Tv_k}.
        \ealigned
      \eequationNN
    \EndIf
    \State \textbf{terminate} and \textbf{return} $(H_{k+1},W_{k+1})$.
  \ealgorithmic
\ealgorithm

\blemma\label{lem.hessian}
  The properties pertaining to $\{H_k\}$ and $\{W_k\}$ stated in Lemma~\ref{lem.well_posed} hold.
\elemma
\proof{Proof.}
  Positive definiteness of $H_{k+1}$ and $W_{k+1}$ follows by induction and the fact that $H_0 \succ 0$ and $W_0 \succ 0$.  In particular, if $s_0 = 0$ or \eqref{eq.Gy_more_d_and_alpha_big} does not hold in iteration $k=0$, then $H_1 \gets H_0 \succ 0$ and $W_1 \gets W_0 \succ 0$; otherwise, positive definiteness of $H_1$ and $W_1$ follows the fact that \eqref{eq.phi} implies $s_0^Tv_0$ and from well-known properties of BFGS updating.  Inductively, positive definiteness of $H_{k+1}$ and $W_{k+1}$ for any $k \in \N{}$ follows by the same arguments.  Finally, with respect to the properties of $\{W_k\}$, the proof follows in a similar manner as that for \cite[Corollary~3.2]{CurtRobiZhou19}.
\endproof

\subsection{Sample point generation}\label{sec.sample}

With $x_{k+1} \in \Dcal$ in hand, a GS algorithm turns to setting the set of sample points $\Xcal_{k+1}$ and corresponding size indicator $p_{k+1}$.  To limit the size of the sample set, which has the benefit of reducing the costs of subsequent QP solves, we follow the lead of \cite{CurtQue15}, which sets $\Xcal_{k+1} \gets \{x_{k+1}\}$ if, for $(\xi,\underline\alpha) \in (0,\infty)^2$, one finds that \eqref{eq.Gy_more_d_and_alpha_big} holds.  Otherwise, the sample set preserves points near $x_{k+1}$ and augments it with randomly generated points, the hallmark of GS methods; see Algorithm~\ref{alg.sample}.

\algdef{SE}[DOWHILE]{DoWhile}{EndDoWhile}{\algorithmicdo}[1]{\algorithmicwhile\ #1}
\balgorithm[ht]
  \caption{Sample Set Update}
  \label{alg.sample}
  \balgorithmic[1]
    \Require $\underline\alpha$ from outer algorithm; $\xi \in (0,\infty)$; $\overline{p} \in \N{}$ with $\overline{p} \geq 1$.
    \If{\eqref{eq.Gy_more_d_and_alpha_big} holds}
      \State set $\Xcal_{k+1} \gets \{x_{k+1}\}$ and $p_{k+1} \gets 0$, \textbf{terminate}, and \textbf{return} $(\Xcal_{k+1},p_{k+1})$.
    \EndIf
    \DoWhile
      \State set $\Scal_{k+1}$ as a set of $\overline{p}$ points from a uniform distribution over $\Bmbb(x_{k+1},\epsilon_{k+1})$.
    \EndDoWhile{$\Scal_{k+1} \not\subset \Dcal$}
    \State Set $\Xcal_{k+1} \gets \{x_{k+1}\} \cup (\Xcal_k \cap \Bmbb(x_{k+1},\epsilon_{k+1})) \cup \Scal_{k+1}$ and $p_{k+1} \gets |\Xcal_{k+1}|-1$.
    \If{$p_{k+1} > p$}
      \State remove the $p_{k+1}-p$ eldest members of $\Xcal_{k+1}$ (except $\{x_{k+1}\}$) and set $p_{k+1} \gets p$.
    \EndIf
    \State \textbf{terminate} and \textbf{return} $(\Xcal_{k+1},p_{k+1})$.
  \ealgorithmic
\ealgorithm

Like for Algorithm~\ref{alg.perturb}, one finds that Algorithm~\ref{alg.sample} can fail if its \textbf{do-while} loop iterates infinitely.  However, under Assumption~\ref{ass.f}, this occurs with probability zero.  The subroutine runs successfully with probability one.

A proof of Lemma~\ref{lem.well_posed}(e) can be found in that for \cite[Lemma~2.5]{CurtQue15}.

\end{document}